\documentclass[english]{article}
\usepackage[T1]{fontenc}
\usepackage[utf8]{inputenc}
\usepackage{morefloats}
\usepackage[letterpaper]{geometry}
\usepackage{float}
\usepackage{amsmath}

\usepackage{graphicx}
\usepackage{color}
\usepackage{adjustbox}
\usepackage{diagbox}
\usepackage{multirow}
\usepackage{array}
\usepackage{mathrsfs}
\usepackage{amsmath, amssymb, amsthm}
\usepackage{subfigure}
\usepackage{subfiles}
\usepackage{tikz}

\usepackage{epstopdf}
\graphicspath{{./figures/}}

\newtheorem{theorem}{Theorem} % we use joint numbering
\newtheorem{definition}[theorem]{Definition}
\newtheorem{assumption}[theorem]{Assumption}

\usepackage{amsfonts} 
\usepackage{color}
\usepackage{adjustbox}
\usepackage{diagbox}
\definecolor{black}{rgb}{0,0,0}

\definecolor{red}{rgb}{1,0,0}

\definecolor{blue}{rgb}{0,0,1}

\newcommand{\cC}{\mathcal{C}}

\newcommand{\cG}{\mathcal{G}}

\newcommand{\cU}{\mathcal{U}}

\newcommand{\cL}{\mathcal{L}}

\newcommand{\cN}{\mathcal{N}}
\newcommand{\cP}{\mathcal{P}}

\newcommand{\cS}{\mathcal{S}}
\newcommand{\bc}{\mathbf{c}}

\newcommand{\bbC}{\mathbb{C}}
\newcommand{\bbR}{\mathbb{R}}

\usepackage{multirow}

\newcommand{\sch}{Schr\"{o}dinger }

\makeatother

\usepackage{babel}
\usepackage{authblk}

\usepackage{graphicx}
\def\al#1\eal{\begin{align}#1\end{align}}
\def\als#1\eals{\begin{align*}#1\end{align*}}
\def\eq#1\eeq{\begin{eqnarray}#1\end{eqnarray}}
\def\eqs#1\eeqs{\begin{eqnarray*}#1\end{eqnarray*}}
\def\eqn#1\eeqn{\begin{equation}#1\end{equation}}
\def\eqns#1\eeqns{\begin{equation*}#1\end{equation*}}
\def\ad#1\ead{\begin{aligned}#1\end{aligned}}
\def\itm#1\eitm{\begin{itemize}#1\end{itemize}}
\def\df#1\edf{\begin{definition}#1\end{definition}}
\def\sl#1\esl{\begin{solution}#1\end{solution}}
\def\ex#1\eex{\begin{example}#1\end{example}}
\DeclareUnicodeCharacter{00A0}{ }
\providecommand{\keywords}[1]
{
	\small	
	\textbf{\textit{Keywords---}} #1
}
\title{}

\title{A clever neural network in solving inverse problems of Schr\"{o}dinger equation}
\author{Yiran Wang\thanks{Department of Mathematics, Purdue University,150 N University St, West Lafayette, IN 47907 }
}

\begin{document}
\maketitle
\begin{abstract}
	In this work, we solve inverse problems of nonlinear Schr\"{o}dinger equations that can be formulated as a learning process of a  special convolutional neural network. Instead of attempting to approximate functions in the inverse problems, we embed a library as a low dimensional manifold  in the network such that unknowns can be reduced to some scalars. The  nonlinear Schr\"{o}dinger equation (NLSE) is $i\frac{\partial \psi}{\partial t}-\beta\frac{\partial^2 \psi}{\partial x^2}+\gamma|\psi|^2\psi+V(x)\psi=0,$ where  the wave function $\psi(x,t)$ is the solution to the forward problem and the potential $V(x)$ is the quantity of interest of the inverse problem. The main contributions of this work come from two aspects. First, we construct a special neural network directly from the Schr\"{o}dinger equation, which is motivated by a splitting method. The physics behind the construction enhances  explainability of the neural network. In particular, each convolution layer  and activation function can be related to corresponding parts of the equation. Under this construction, it can be rigorously proved that the neural network has a  convergence rate with respect to the length of input data and number of layers. The other part is using a library-search algorithm to project the solution space of the inverse problem to a lower-dimensional space.  The way to seek the solution in a reduced approximation space can be traced back to the compressed sensing theory.  The motivation of this part is to alleviate the training burden in estimating functions. Instead, with a well-chosen library, one can greatly simplify the training process.  More specifically, in one of the experiments, we analysis the landscape of the loss function with respect to the training parameters to help us  obtain the optimal solution to the inverse problem. A brief analysis  is given, which focuses on  well-possedness of some mentioned inverse problems and convergence of the neural network approximation. To show the effectiveness of the proposed method, we explore in some representative problems including simple equations and a couple equation. The results can well verify the theory part. In the future, we can further explore manifold learning to enhance the approximation effect of the library-search algorithm.
	\end{abstract}

\keywords{Physics-inspired neural network, Reduced order modeling, Inverse problems, Schr\"{o}dinger equation, Compressed sensing}
\section{Introduction}
Neural network approach has been widely explored in solving partial differential equations (PDE) since the pioneer work \cite{lagaris1998artificial}. Most of them choose from some existing architectures and train the model to obtain satisfying results for the PDEs. The advantage is that it does not require too much information from the equation itself especially when it is a  data-driving method. This provides great convenience in solving different types of equations hence it is user-friendly, which also brings some drawback on the other hand. Performing as a blackbox, this kind of neural network approach lacks sufficient explainability, which comes from three aspects. The first one is the existence of parameters that are optimal to some sense from theoretical aspect.  This is somehow answered by the popular universal approximation theorem for neural networks. However, it is not sufficient to choose only one single layer in many cases. For a network with a particular architecture,  the existence of optimal parameters is unclear. Secondly, even though one can obtain the optimal solution from theory, it remains unclear whether the optimal solutions can be attained with the training process. In many cases, one can only find it with trial and error. The third uncertainty exists in the generalization. If one can obtain a set of ``good'' parameters to build an optimal neural network, whether the performance is satisfying when different dataset are used. If it performs well in some similar data as the training data,  which function space it can approximate.  Overall, there is not a general explanation about the above three  problems, i.e., the existence of optimal parameters, how to train it, and the generalization power. 

To enhance the explainability, one can equip the construction of neural network with  information of equations instead of treating networks and PDEs separately.  To this end, physics-informed neural newtork (PINN) \cite{raissi2019physics} serves as a representative example by utilizing a loss function derived from equation. Inspired by this work, we  reconsider the relationship of neural network and PDEs and desire to figure out a novel structure of network that makes use of the equation itself.  As for the interested equations, we  focus on a class of nonlinear Schr\"{o}dinger equation since they are significant in both scientific research and industry \cite{guo2006some,neuhasuer1989time,kohn1954solution,scalora2005generalized}. Moreover, numerical treatments of Schr\"{o}dinger equation have been  investigated extensively \cite{bandrauk1993exponential,chang1999difference,stocker1999current,mills2017deep,bao2003explicit}, among which time splitting schemes are widely  used  \cite{bao2002time,bao2003numerical,thalhammer2012convergence,besse2002order,thalhammer2008high}.   Apart from forward process, inverse problems attract attention from many researchers due to its wide applications \cite{mueller2012linear} and contribute some effective numerical techniques \cite{mandache2001exponential,backus1967numerical,razavy2020introduction,ablowitz2007inverse,snieder1998role}. Based on the understanding of  the time splitting method,  we consider designing a novel network from the nonlinear Schr\"{o}dinger equation.  Hence, we call the derived network  NLS-Net. In \cite{wang2021inverse}, we inverstigated an coefficient inverse problem in nonlinear Schr\"{o}dinger equation without potential. In this work, we explore an inverse problem of finding potential from data with the NLS-Net. 

Compared with scalars, it is far more complicated to seek functions. To this end, we introduce a library-search approach to reduce the problem to a coefficient inverse problem. This method is motivated from an algorithm in  sparse identification of nonlinear dynamics (SINDy-PI) \cite{kaheman2020sindy}. The priciple idea is to choose some functions from  a reference set, which can be understood as a library. We further use a linear combination of the chosen functions to represent the unknown potential to the inverse problem. In other words, we project the original solution to a fixed approximation space, hence we only need to determine the corresponding coefficients. However, from the perspective of expressive power, we need to incoporate sufficient modes in the library, which performs similarly as utilizing an either deep or wide neural network. Put differently, it is obvious that approximation power comes from sufficient degree of freedom. On the other hand,  increasing the number of unknowns means  a growth in computation cost, which is apparently a side effect. For example, everyone in the area of scientific computing is familiar with the word ``the curse of dimensionality'' \cite{koppen2000curse}, which may bring a blow up in computation cost. Hence, model reduction methods \cite{antoulas2000survey,benner2015survey,penzl2006algorithms,beck1996model,obinata2012model,phillips2004poor} are introduced to seek a ``cheaper'' model based on the original model to attain the quantity of interest.  In particular,  there are considerable number of works in reduced order methods of solving PDEs numerically \cite{nouy2010priori,benner2017model,xiao2017parameterized,audouze2013nonintrusive,nguyen2008efficient,volkwein2013proper,ravindran2000reduced,gunzburger2007reduced}. For inverse problems in PDEs, many researchers devote good ideas of applying model reuction \cite{frangos2010surrogate,cui2015data,galbally2010non,lieberman2010parameter,fu2018pod,ghattas2021learning}. It is understandable that one needs to sacrifice accuracy in order to reduce the computation cost since there is no free lunch. However if one has prior information about the quantity of interest, it is possible to use some well-designed models that can attain a relatively high accuracy with reduced workload. In image or signal processing, compressed sensing \cite{miosso2009compressive,donoho2006compressed,tsaig2006extensions,duarte2011structured,lustig2008compressed,shi2009advances,razzaque2014energy} is utilized for enforcing sparsity in some incorporated variables based on the knowledge that the solution lies in a low-dimensional subpace. Therefore, it is reasonable to apply this technique in the area of numerical PDE as a other reduced order approach, which becomes another motivation of this work where we desire to  approximate the unknown potential in the nonlinear Schr\"{o}dinger equation. The prior information is that the potential function should have a relatively simple expression, which can be a linear combination of a few different basis functions. Hence,  it is reasonable to assume the target coefficient vectors with respect to bases exploits a sparse form, which explains the motivation of applying compressed sensing in this context. Overall,  contributed by above motivations, we propose  a library-search approach combined with compressed sensing to train a special network called NLS-Net that is derived from NLSE. The proposed method is potential to be developed a new framework of solving inverse problems. To evaluate the performance of the proposed method, we use three reprentative examples, two of which are single equations with unknown potentials while the rest one is coupled system of two NLSE. In all the examples, we use NLS-Net to obtain the forward solutions.  From  \cite{bao2012mathematical},  one can derive an upper bound for the approximation errors of forward solutions, which can in turn guarantee the convergence of NLS-Net with respect to number of layers and data amount. In our numerical results, we show that the performance of the NLS-Net is consistent to the theory. For the inverse part, we use the library-search combined with compressed sensing to seek unknown potentials. In view of the specificity of the third case, we apply a landscape analysis to help us find the optimal solution more efficiently. It is satisfying that the concerned loss as well as errors present an exponential decay in all the three examples. Moreover, we compare the numerical solutions of the inverse part, which are approximately identical to the groundtruth. Both of them shows the effect and accuracy of the proposed method.

The paper is organized in six parts. We present the concerned problem in the first part and review the conventional split-step methods in the next part. The proposed method is introduced in the third section and it is compsed of four subsections. We first introduce the mechanism of designing a convolutional neural network from NLSE. Following this, we provide some background in compressed sensing. The library-search algorithm is then proposed, which is combined with compressed sensing to form the final approach in the last part. We alsp present some relevant theroy in the fifth section. Numerical simulations are demonstated in the following part. We conclude this work and offer some remark and potential interested problems in the last place.

\section{Problem description}
From an initial value  problem of Schr\"{o}dinger equation as follows:
\begin{equation}
	\begin{aligned}
i\frac{\partial \psi}{\partial t}+\beta\frac{\partial^2 \psi}{\partial x^2}+\gamma|\psi|^2\psi+V(x)\psi&=0,\\
\psi(0,x)&=h(x),\label{model_ori}
\end{aligned}
\end{equation}
where $\psi(x,t): \bbR\times [0,T]\to \bbC$ is the wave function and $V(x)$ is the potential to be estimated. Moreover, $\beta$,$\gamma$ are constants and $h(x):\bbR\to \bbC$ is a given as the initial condition. The observed dataset composed of continuous functions is defined to be $\mathcal{D}_{\text{cont}}:=\{\psi_0(x),\psi_T(x)\}$ where $\psi_0(x)$ and $\psi_T(x)$ is the observed wave at the initial time $0$ and final time $T$, respectively. Based on the observed data, we estimate the potential $V(x)$. In particular, we define the $S(t;V)$ to be \sch operator in \eqref{model_ori} as $\psi(\cdot,t)=S(t;V)\psi(0,\cdot)$. To identify the potential $V$ from $\mathcal{D}_{\text{cont}}$, we formulate the following optimization problem:
\al
\min\limits_{V}\|S(t;V)\psi_0-\psi_T\|_2^2. \label{opt_ori}
\eal
Numerically, we can not use the exact $S(t;V)$ operator that requires directly solving \eqref{model_ori} which is infeasible. Instead, we use operator splitting \cite{holden2010splitting,leveque1982time,holden2013operator} as an approximation method.
\section{Conventional split-step methods}
In this part, we review some theoretical results on the well-posedness and stability of the inverse problem. Besides, we state the conventional split-step method.
First, we rewrite \eqref{model_ori} into the form
\begin{align}
	\frac{\partial}{\partial t}\psi(x,t) = (\mathcal{L}+\mathcal{N})\psi(x,t),\quad x\in \bbR.\label{model_split}
	\end{align}
Here $\cL:=-i\beta\frac{\partial}{\partial x^2}$ is the linear operator and $\cN:=i\gamma|\psi(x,t)|^2+iV(x,t)$ is the nonlinear operator. Instead of directly solving \eqref{model_split}, one can split it into two subproblems corresponding to two operators $\cL$ and $\cN$, respectively. In particular, one can  solve a linear problem 
\al
\frac{\partial \psi}{\partial t} = \cL\psi_L(x,t), \label{sub_pro1}
\eal
and a nonlinear problem
\al
\frac{\partial \psi}{\partial t} = \cN\psi_N(x,t).\label{sub_pro2}
\eal
Both of the above two problems can be solved sequentially to obtain a final solution. We denote $\delta t$  to be the time step size . In terms of performing the splitting, there are different orders of the method. The simplest version is the first-order method, in which is each advancement of $\Delta$ is composed of two steps. First, one solves either the linear \eqref{sub_pro1} or the nonlinear problem \eqref{sub_pro2} to generate a intermediate solution. Then one apply the other operator to finish this time block. More precisely, if one choose $\cL$ prior to $\cN$ for each $\Delta t$, then a direct result after forwarding one time block is 
\als
 \psi(x,t+\Delta t) \approx \exp(\cN\Delta t)\exp(\cL\Delta t)\psi(x,t).
\eals
Similarly, the second-order method also be decoupled into alternatively solving linear \eqref{sub_pro1} and \eqref{sub_pro2}. The only difference is we split each $\Delta$ into three steps. The key idea is that the solution of one step serves as the initial condition for the next step. The most popular second-order splitting is the standard Strang splitting with the form 
\als
 \psi(x,t+\Delta t) \approx\exp(\cL\Delta t/2) \exp(\cN\Delta t)\exp(\cL\Delta t/2)\psi(x,t).
\eals
The advantage of the splitting method is that both \eqref{sub_pro1} and \eqref{sub_pro2} can be much easier than \eqref{model_split}. In particular, one can explicitly solve \eqref{sub_pro2} and for \eqref{sub_pro1}, many choices can be chosen, for example, finite difference schemes or fourier method. The final combination of splitting and spatial discretization ways form the different numerical methods like split-step finite difference (SSFD) and split-step fourier method (SSFM).

\section{A physics-inspired neural network with a library-search approach }
Compared the last section, we now introduce an new approach with a special-designed neural network (NN) \cite{wang2021inverse} and a library-search way to solve some interested inverse problems. First of all, we first give an introduction of the special NN. Then, we introduce the library-search approach which serves as a low-rank representation of the solution to the inverse problem. Finally, we summarize the overall method in the last subsection.
\subsection{Learning process of a convolutional neural network}
Motivated by SSFM, we can formulate a convolutional neural network to solve the inverse problem. In particular, for linear and nonlinear operator in \eqref{sub_pro1} and \eqref{sub_pro2}, we use a convolutional operator and a nonlinear operator, whose can be explicitly derived from the equation. As the solution of the inverse problem is incorporated into the NN layers, the backpropagation essentially serves as a process of the seeking the ground truth of the inverse problem. We define the linear and nonlinear part as follows. 

For linear problem, we apply 
\begin{align}
\psi_L(x,t+\Delta t) = \kappa \int K(s;\beta\Delta t) \psi_L(x-s,t),\label{cont_lin_op}
\end{align}
where $K(s;\eta)$ is a convolution kernel defined by
$$K(s;\eta):= \sqrt{\frac{i}{\pi \eta}}\exp\left(-\frac{is^2}{4\eta}\right), \forall s\in \bbR.$$ 
 Furthermore, we can introduce a convolution operator $\cC(\cdot;\beta\Delta t)$ such that 
 \begin{align}
 	\psi_L(\cdot,t+\Delta t) = \mathcal{C}(\cdot;\beta\Delta t)	\psi(\cdot,t):=  K(s;\beta\Delta t) *\psi_L(\cdot,t),\label{cont_lin_op2},
 \end{align}
where $*$ refers to convolution with respect to $x$.
For nonlinear problem, we use 
\begin{align}
\psi_N(x,t+\Delta t) = \cG(\psi_N;\gamma\Delta t,V\Delta t),\label{cont_nonl_op}
\end{align}
where the nonlinear operator $\cG$ is defined as
$$\cG(u;g,h):=u\exp{i(g |u|^2+h)}.$$
Finally we formulate the forward propagation as 
\begin{align}
	\psi(\cdot,t+N\Delta t)\approx \cS(\cdot;\beta,\gamma,V,\Delta t)\psi(\cdot,t):=
\left(\cC\left(\cdot;\frac{\beta\Delta t}{2}\right) \circ \cG(\cdot;\gamma\Delta t,V\Delta t) \circ \cC\left(\cdot;\frac{\beta\Delta t}{2}\right)\right)^{\circ N}\psi(\cdot,t)\label{whole_cont}
	\end{align}
Corresponding to the continuous version of linear and nonlinear operators in \eqref{cont_lin_op} and \eqref{cont_nonl_op}, we further define the discrete version of $\cL$ and $\cN$.  Before this, we need to consider a corresponding equtaion with periodic boundary condition as follows. 
\begin{equation}
	\begin{aligned}
		i\frac{\partial \psi}{\partial t}+\beta\frac{\partial^2 \psi}{\partial x^2}+\gamma|\psi|^2\psi+V(x)\psi&=0, \quad x\in \Omega, t\in [0,T],\\
		\psi(x,0)&=h(x),\quad x\in \Omega,\\
		\psi(a,t)&=\psi(b,t),\quad t\in [0,T]\\
		\frac{\partial \psi}{\partial x}(a,t) &= \frac{\partial \psi }{\partial x}(b,t),\quad t\in [0,T],
		\label{model_periodic_bc}
	\end{aligned}
\end{equation}
where $\Omega$ is specified as $[a,b]$.
We choose homogeneous grid $\{a_j\}_{j=1}^{M}$ with  $a_j = a+\frac{j(b-a)}{M}$ for $j=1,\ldots, M.$ We define $\Delta x = \frac{b-a}{M}$. Since we enforce the periodic boundary condition, i.e., $\psi(t,a)=\psi(t,b)$, we choose $b$ in the grid and $a$ is not in the grid.
Then, we define a discrete convolution kernel as 
\begin{align*}
	\tilde{K}(q;\eta)[j]:=\sqrt{\frac{i}{\pi \eta}}\exp\left(-\frac{i(jq)^2}{4\eta}\right), \quad \forall j\in [M],
	\end{align*} 
where $M$ is the number of collocation points in the spatial domain and we use a combinatorial notation $[M]={1,\ldots,M}$. We define $t_n = n\Delta t$. With the spatial discretization into the grid $\{a_j\}_{j=1}^{M}$ , we define $\tilde{\psi}_{M,L}^m$ to be the approximate solution to \eqref{cont_lin_op}  at time $t_{m}$, from which we update the solution at $t_{m+1}$ with discrete convolution as follows,
\begin{align*}
	\psi_L(x_n,t_{m+1})\approx \tilde{\psi}_{M,L}^{m+1}[n]=\tilde{K}(q;\beta\Delta t)*\tilde{\psi}_{M,L}^m:=\sum_{j=1}^{M}\tilde{K}(q;\beta\Delta t)[j]\tilde{\psi}_{M,L}^m[n-j],
\end{align*}
where $\tilde{K}(q;\eta)$ are replaced by zero when $n-j\notin [M]$. Again, we define an associated discrete convolution operator  $\tilde{\cC}(\cdot ;\eta)u:=\tilde{K}(q;\eta)*u$. As for the nonlinear operator, we use $\tilde {\cG}$ for the discrete version but the definition of $\tilde {\cG}$ is the same as $\cG$ but we can restrict the domain to set of discrete vector.
Corresponding to \eqref{whole_cont}, we can summarize the discrete version by combining the $\cC$ and $\cG$ as 
\begin{align}
\tilde{\psi}_{M}^{N} = \tilde{\cS}(\cdot;\beta,\gamma,V,\Delta t)\tilde{\psi}_M^{0}:=\left(\tilde{\cC}\left(\cdot;\frac{\beta\Delta t}{2}\right) \circ \tilde{\cG}(\cdot;\gamma\Delta t,V\Delta t) \circ \tilde{\cC}\left(\cdot;\frac{\beta\Delta t}{2}\right)\right)^{\circ N}\tilde{\psi}_M^{0},\label{whole_disc}
	\end{align}
where we define $\tilde{\psi}_{M}^{N}$ as the discrete solution to \eqref{model_periodic_bc} with the operator splitting method. 
This can be expressed as a $N-$ layer complex-valued CNN with activation function $\tilde{\cG}$ as below:
\begin{align}
	{\tilde{\mathcal{W}}}(\Theta)=\widetilde{\mathcal{C}}_h \circ \cG \circ \widetilde{\cC} \circ \ldots \circ \widetilde{\cC} \circ  \cG \circ \widetilde{\cC}_h. \label{NN}
	\end{align}
Here $\Theta$ is the set of weights contained in the operators $\tilde{\cC}$ and $\tilde{\cG}$ and $\widetilde{\mathcal{C}}_h$ refers to $\frac{\Delta t}{2}$ of forwarding in the linear step.
\begin{figure}[htbp!]
	\centering
	\includegraphics[width=1.0\textwidth]{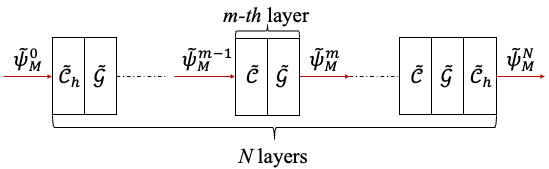}
	\caption{An illustration of NLS-Net composed of $N$ layers and each layer is a composition of $\tilde{\cC}$ and $\tilde{\cG}$. Here we use $\tilde{\cC}_h$ denotes $\tilde{\cC}\left(\cdot;\frac{\beta\Delta t}{2}\right)$; $\tilde{\cC}$ and  $\tilde{\cG}$ refers to $\tilde{\cC}(\cdot;\beta\Delta t)$ and $\tilde{\cG}(\cdot;\gamma\Delta t,V\Delta t)$.}
	\end{figure}

\subsection{Compressed sensing  background}
In \cite{donoho2006compressed}, the problem is to reconstruct unknown image or signal $x\in \bbR^m$  from $n$ measurements. Based on the prior information that $x$ has a sparse representation in some orthonormal basis, the number of measurements $n$ can be dramatically smaller than the size $m$.  For example, some natural classes of images with $m$ pixels need only $n = O(m^{1/4}\log^{\frac{5}{2}}(m))$ samples for relatively accurate recovery. Here, we briefly state the problem and focus on the information that can be used in the proposed method.

Suppose $X$ is a class of objects of interest and one needs to recover $X$ from some measurements. First  they define an information operator $I_n: X\to \mathbb{R}^n$ as follows:
\begin{align*}
	I_n(x)=(\langle \xi_1,x\rangle,\ldots, \langle \xi_n,x\rangle)
	\end{align*}
where $\xi_i$ are sampling kernels. Furthermore, $A_n:\mathbb{R}^n\to \mathbb{R}^m$ is an algorithm that reconstruct $x$ from $I_n(x)$. Overall, the problem is formulated as minimizing the $l_2$ error of the reconstruction
\begin{align*}
E_n(X) = \inf\limits_{A_n, I_n} \sup\limits_{x\in X}\|x-A_n(I_n(x))\|_2.
\end{align*}
Here, $A_n$, $I_n $ can be any possible methods of linear sampling and reconstruction. The goal is to obtain a faithful recovery with limited measurements.   In particular, it can be splite into two questions: how to choose $I_n$ and once it is fixed, how to determine an effective way to reconstruct $x$, i.e., finding $A_n$.  First of all, the following assumption is used.
\begin{assumption}
	For each $x\in X$, there exists a set of orthogonal bases $(\psi_i\in \bbR^m, i=1,\ldots, m)$. Define $\theta_i := \langle x, \psi_i\rangle$ to be the corresponding transform coefficients as $\theta_i := \langle x, \psi_i\rangle$. Then $x$ is sparse in the sense that  for some $0<p<2$ that 
	\begin{align*}
		\|\theta\|_p\equiv \left(\sum\limits_{i} |\theta_i|^p\right)^{\frac{1}{p}}\leq R
	\end{align*}
	for some $R>0$.
	\end{assumption}
The authors are only interested in near-optimal performance, i.e., optimal up to some constants. In \cite{donoho2006compressed}, they provided the condition of near-optimal information operator and the reconstruction procedure based on linear programming. Briefly speaking, $\Phi$ is an $n\times m$ matrix generated by randomly sampling the columns and $J\subset \{1,\ldots, m\}$ is a set of index. Let $\Phi_j$ denote a submatrix of $\Phi$, which is a result of choosing indicated columns of $\Phi$. This matrix contains the information of sampling.   We pay special attention to one condition of $\Phi_j$ for the desired approximation power: for strict positive parameter $\rho$, the minimal  singular value of $\Phi_J$ exceeds $\eta_1>0$ uniformly in $|J|<\frac{\rho n}{\log(m)}$. Here, the columns of $\Phi_j$ are used as bases functions. To attain an effective approximation, it is significant to obtain a particular quantitative degree of linear independence among bases. Based on the information matrix $I_n$, one can solve the following convex optimization problem:
\begin{align*}
	\min\|\Psi^T\|_1 \text{ subject to } y_n = \Phi\Psi^T(x),
	\end{align*}
where $\Psi\in \mathbb{R}^{m\times n}$ is an orthorgonal matrix. In other words, a desired form is to find $x$ with smallest $l_1$ norm of transform coefficients and consistent with the measurements $y_n$.

\subsection{A library-search algorithm}
Motivated by the \cite{kaheman2020sindy}, we can formulate a library-search algorithm in solving the inverse problem. Suppose we have a set of functions denoted by $\cP = \{g_1,\ldots,g_{N_L}\}$, where $N_L$ is the cardinality of the function set $\Phi$. The key idea of this library-search method is that we solve $V$ in a space spanned by $\cS$. Mathematically, we can formulate the process as follows. Given a set of collocation points $\{x_1,\ldots,x_{N_x}\}$, we define vectors $g_{i,d}\in \bbR^{N_x}$ by $g_{i,d}[j]:=g_i(x_j)$, $\forall j\in [N_x]$, $\forall i\in [N_L]$. The discretized target solution $\psi_T$ can be similarly defined as $\psi_{T,d}[j] = \psi_{T}[x_j]$, $\forall j \in [N_x]$. Now we have obtained a set of bases $\{g_{i,d}\}$, which can generate a subspace $H$ in $\bbR^{N_x}$.  We further assemble a matrix $\Phi \in \bbR^{N_x\times N_L}$ define by $\Phi e_j = g_{i,d}$, where $e_j$ are the $j$-th standard basis for $\bbR^{N_L}$. Given a coordinate vector $\bc\in \bbR^{N_L}$, one can generate a vector $y\in H$ defined by $y=\Phi \bc$.
Now we replace \eqref{opt_ori} by the following problem:
\al
\min\limits_{\bc}\|\tilde{\cS}(\cdot;\beta,\gamma,\Phi \bc,\Delta t)\tilde{\psi}_L^{0}-\psi_{T,d}\|_2^2. \label{opt_new}
\eal
\subsection{Our approach}
From first subsection, we explore some details of the construction of the convolutional neural network derived from the Schr\"{o}dinger equation, which is therefore denoted by NLS-Net. To represent the potential, we introduce the library-search algorithm, which indeed serves as an effective way of parameterizing the NLS-Net. Inspired by the compressed sensing, we can seek an approximation with relatively high sparsity. With the notation in the compressed sensing, we desire to find a way to find a corresponding transform coefficients with relatively small $l_1$ norm.

For the training process of NLS-Net, we use a summation of $l^2$ error and $l^1$ norm of $\bc$ for the loss function defined below:
\al
J_{\bc} =\displaystyle\frac{1}{2} \frac{\|\tilde{\cS}(\cdot;\beta,\gamma,\Phi \bc,\Delta t)\tilde{\psi}_L^{0}-\psi_{T,d}\|_2^2}{\|\psi_{T,d}\|_2^2}+\lambda\|c\|_1. \label{loss}
\eal
To minimize $J_{\bc}$, we use proximal gradient descent \cite{nitanda2014stochastic}. We summarize the method of combining the NLS-Net with the library search algorithm in the following Table.
\begin{table}
	\centering
	\begin{tabular}{c|l}
		\hline \hline
		\textbf{Input}&$\bc_{\text{ini}}$, $\Phi$, $\beta$, $\gamma$, $\Delta t$,$\tau$\\
		\textbf{Initialization:}& Build the NLS-Net and initializing parameters, where $\bc = \bc_{\text{ini}}$.\\
		\textbf{Step 1:}& Compute the loss $J_{\bc}$ defined in \eqref{loss}.\\
		\textbf{Step 2}& If $J_{\bc}\leq \tau$, output result. Else, continue step 3 and 4.\\
		\textbf{Step 3:}& Update $\bc$ based on backpropagation and return to Step 1.\\
		\textbf{Output}&$V=\Phi \bc$\\
		\hline
	\end{tabular}
\end{table}

\section{Theoretical background}
In this section, we present the theoretical background of the inverse problem. We briefly mention two parts as follows:
\begin{enumerate}
	\item The well-posedness of the inverse problem.
	\item The convergence analysis of the NLS-Net
\end{enumerate}
In particular, in the first part, we review some important results about uniqueness and stability of a special form of inverse problem. We remark that the existing result is limited to linear case however, it offers insights of the nonlinear case which deserves exploration. For the second part, we present the convergence rate of the constructed neural network with respect to the number of collocation points and layers. \\
\\
In \cite{baudouin2002uniqueness}, Baudouin et al. (2007) proved uniqueness and stability in an inverse problem for the Schr\"{o}dinger equation. They consider the following equation:
\al
\left\{\begin{array}{l}
	i y^{\prime}(x, t)+\Delta y(x, t)+q(x) y(x, t)=0, x \in \Omega, t \in(0, T) \\
	y(x, t)=h(x, t), x \in \partial \Omega, t \in(0, T) \\
	y(x, 0)=y_0(x), x \in \Omega .
\end{array}\right.\label{paper1}
\eal
For this equation, they proved the following theorem.
\begin{theorem}
	Let $\cU$ be a bounded subset of $L^{\infty}(\Omega)$, $q\in L^{\infty}(\Omega)$ and $y$ be a  solution of \eqref{paper1}. Assuming $$
	\begin{gathered}
		\exists x_0 \in \mathbb{R}^N \backslash \bar{\Omega} \text { such that } \Gamma_0 \supset\left\{x \in \partial \Omega ;\left(x-x_0\right) \cdot \nu(x) \geq 0\right\} \\
		y(q) \in W^{1,2}\left(0, T, L^{\infty}(\Omega)\right)
	\end{gathered}
	$$
	$y_0$ is real valued and $\left|y_0\right| \geq r_0>0$, ae in $\bar{\Omega}$.
	There exists a constant $C=C(\Omega,T,\Gamma_0,\|q\|_{L^{\infty}(\Omega)},y_0,h,\cU)>0$ such that if 
	$$\frac{\partial y(q)}{\partial \nu}-\frac{\partial y(p)}{\partial \nu}\in H^1(0,T;L^2(\Gamma_0))$$
	then $\forall p\in\cU,$
	\als
	\|q-p\|_{L^2(\Omega)}\leq C\left\|\frac{\partial y(q)}{\partial \nu}-\frac{\partial y(p)}{\partial \nu}\right\|_{H^1(0,T;L^2(\Gamma_0))}.
	\eals
\end{theorem}
This theorem characterizes the fact that the solutions to the inverse problem change continuously with the normal derivative of the forward solutions.

In the following part, we review some important results of the convergence rate of the split-step method, which is the core of the NLS-Net. In \cite{lubich2008splitting}, Lubich (2008) proved the convergence rate of the splitting methods for Schr\"{o}dinger-Poission and cubic nonlinear Schr\"{o}dinger equations. Here we only consider the latter one as follows:
\al
i\frac{\partial \psi}{\partial t} = -\nabla \psi +V\psi, \quad x\in \bbR^3,\quad t>0, \label{model_SSFM}
\eal
where $V = V[\psi]=\pm|\psi|^2.$  Assuming the solution $\psi(t)$ to the \eqref{model_SSFM} is in $H^4$ and let $m = \max_{0\leq t\leq T}\|\psi(t)\|_{H^4}$, the following theorem is proved.
\begin{theorem}
	The numerical solution $\psi_n$ given by the strang-split scheme with step size $\tau$ for equation \eqref{model_SSFM} has a first-order error bound in $H^2$ and a second-order error bound in $L^2$, for $t_n = n\tau\leq T$,
	\als
	\|\psi_n-\psi(t_n)\|_{H^2}\leq C(m,T)\tau,\\
	\|\psi_n-\psi(t_n)\|_{L^2}\leq C(m,T)\tau^2.
	\eals
\end{theorem}
This theorem provides the convergence rate of strang-split scheme with respect to the time step $\tau$ under $H^2$ and $L^2$ norm. However, it is a continuous case in the spatial domain and they didn't consider a discretized case. Hence, we proceed to present the result given in  \cite{bao2012mathematical} as follows.
\begin{theorem}
 Let $\psi^n \in X_M$ be the numerical approximation obtained by the TSSP \cite{bao2012mathematical}. Under assumptions (4.84) and (4.85), there exist constants $0<\tau_0, h_0 \leq 1$, such that if $0<h \leq h_0, 0<\tau \leq \tau_0$ and $m \geq 5$, we have
	$$
	\begin{aligned}
		& \left\|\psi\left(x, t_n\right)-I_M\left(\psi^n\right)(x)\right\|_{L^2(U)} \lesssim h^m+\tau^2, \quad\left\|\psi^n\right\|_{\infty} \leq M_1+1, \\
		& \left\|\nabla\left(\psi\left(x, t_n\right)-I_M\left(\psi^n\right)(x)\right)\right\|_{L^2(U)} \lesssim h^{m-1}+\tau^2, \quad 0 \leq n \leq \frac{T}{\tau},
	\end{aligned}
	$$
	where the interpolation operator $I_M$ is given below and $M_1=\max _{t \in[0, T]}\|\psi(\cdot, t)\|_{L^{\infty}}$.
\end{theorem}
The mentioned notations are defined as follows.
\begin{enumerate}
	\item 
	$X_M=\{v = (v_j)_{j\in \mathcal{T}_M^0}|v_0=v_M=0\}\subset C^{M+1}$ with $\mathcal{T}_M^0 = \{j| j=0,1,\ldots, M\}$.
	\item 
	TSSP refers to time-splitting sine pseudospectral method introduced in \cite{bao2012mathematical}. The difference of this method from split step fourier method is that the spatial discretization is sine spectral method while the latter one is fourier method. 
	\item 
	$I_M$ is an interpolation operator defined as 
	 \begin{align*}
	 	(I_M\phi)(x) = \sum\limits_{l=1}^{M-1}\tilde{\phi}_l\sin(\mu_l(x-a)),
	 	\end{align*}
 where $\tilde{\phi}_l=\frac{2}{M}\sum\limits_{j=1}^{M-1}\phi_j\sin(\frac{jl\pi}{M})$ with $l\in \mathcal{T}_M$ and $\phi_j = \phi(x_j)$.
	\end{enumerate}
Motivated from  this theorem, we can have the following upper bound for NLS-Net.
\begin{theorem}\label{thm: conv_nn}
	We define that $\tilde{\psi}_M^n$ is the solution at time $t_n = n\Delta t$ generated by \eqref{whole_disc} with $M$ as the input data length and $N$ is the number of layers. We recall that $\psi(x,t_n)$ is the analytical solution to \eqref{model_ori} at time $t_n$. Suppose that the initial value $\psi(\cdot,0)\in H^s(\bbR\to \bbC)$ for $s\geq 5$. Then there exist parameters $\Theta$ of the neural network $\tilde{\mathcal{W}}(\Theta)$, which is defined by \eqref{NN} such that 
	\begin{align}
		\|\tilde{\mathcal{W}}(\Theta)\tilde{\psi}_M^n-\psi_{T,d}\|\lesssim (\Delta x)^s+(\Delta t)^2,
	\end{align}
	where $\Delta x = \frac{b-a}{M}$ and $\Delta t = \frac{T}{N}$.
\end{theorem}
We remark that this theorem specifies that the approximation effect of the neural network converges respect to the spatial and time step sizes. 
\section{Numerical experiments}
In this section, we demonstrate the performance of the proposed method with three examples. In the first example, we use the NLS-Net to approximate the one-soliton solution for the generalized nonlinear Schr\"{o}dinger equation (NLSE), where the potential is exponential. For the second one, we further consider another form of NLSE with a periodic potential, which is called Gross-Pitaevskii (GP) equation. In both examples, only a single equation is involoved. For the third one, we explore a coupled system with two equations, where the potential is a function of solutions to both equations. In all the examples, we transform the problem of function estimating to seeking some coefficients, which are further used to recover the potentials. In particular, in the first examples with single equations, we use the library-search algorithm to construct a low-dimensional approxiamtion of the true potential.  More specifically, we use a library  as follows
$$\Phi:=\{x, \sin(x),\cos(x),x^2,\sin^2(x), \cos^2(x), e^{-x}, e^{-2x}, e^{-x^2},e^{-2x^2}\}.$$ In addition, we make use of compressed sensing to obtain a relatively sparse coefficient vector, with which only very few functions are involved in the final numerical solution. For the third one, owing to the special setting of the problem, one can easily relate the unknown function to some scalars hence the problem can well simplified.  To present the performance of the proposed method, we define $e_{\phi}$ and $e_{V}$ as below:
\begin{align}
	e_{\psi} = \frac{1}{2}\frac{\|\psi^{\text{numerical}}-\psi^\text{exact}\|_2^2}{\|\psi^\text{exact}\|_2^2},\quad   e_{V}=  \frac{\|V^{\text{numerical}}- V^{\text{exact}}\|_2}{\|V^{\text{exact}}\|_2},
	\end{align}
where $\psi^{\text{numerical}}$ and $V^{\text{numerical}}$ are the numerical results where $\psi^{\text{exact}}$ and $V^{\text{exact}}$ are corresponding references. In particular, $e_{\psi}$ is from the  $l^2$ term in loss function $J_{\bc}$, which we show more interest in compared with the whole loss $J_{\bc}$.\\
\newline
\textbf{Example 1.} One-soliton solution for NLSE.

\begin{figure}[htbp!]
	\centering
	\subfigure{\includegraphics[width=0.45\textwidth]{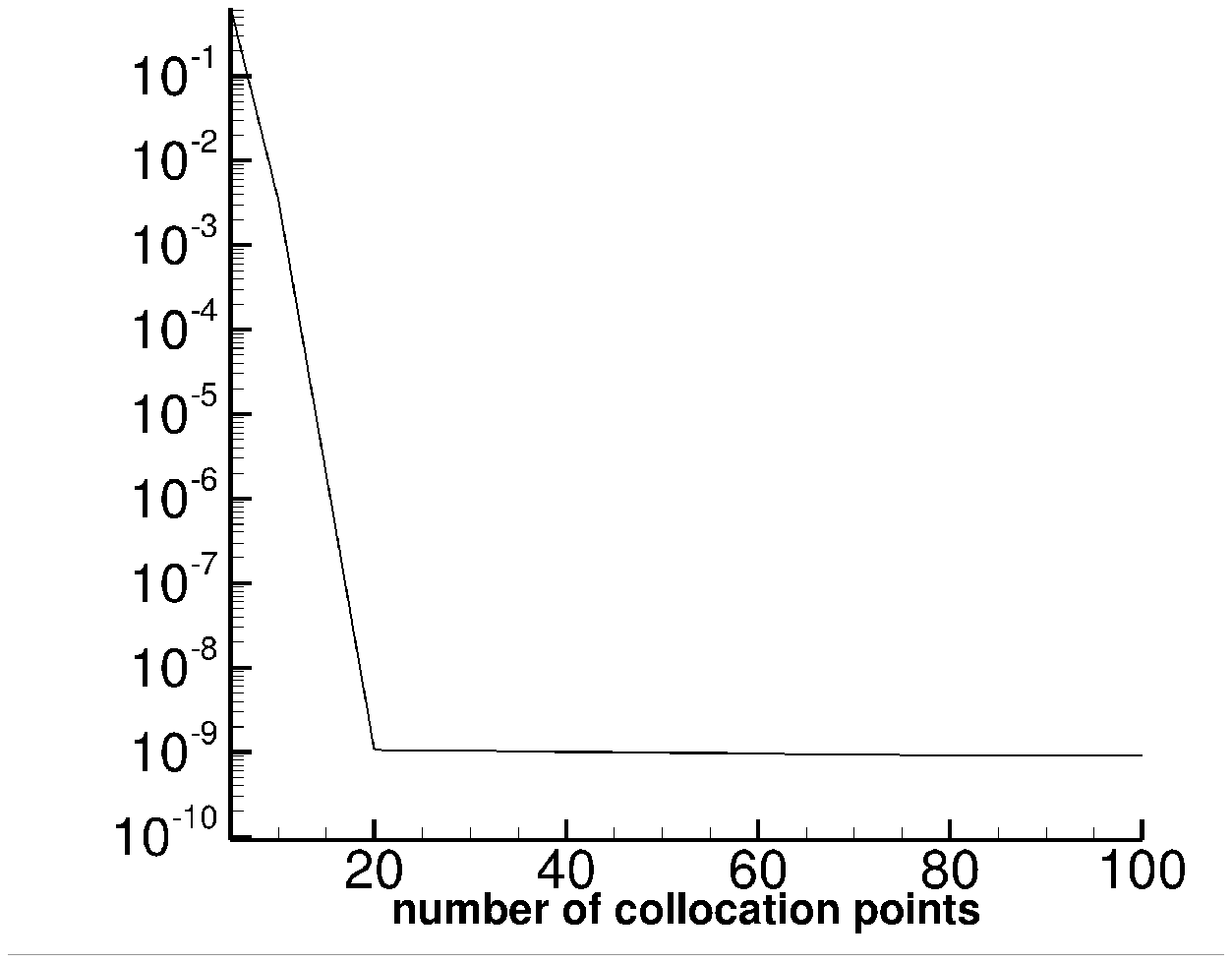}}
	\subfigure{\includegraphics[width=0.45\textwidth]{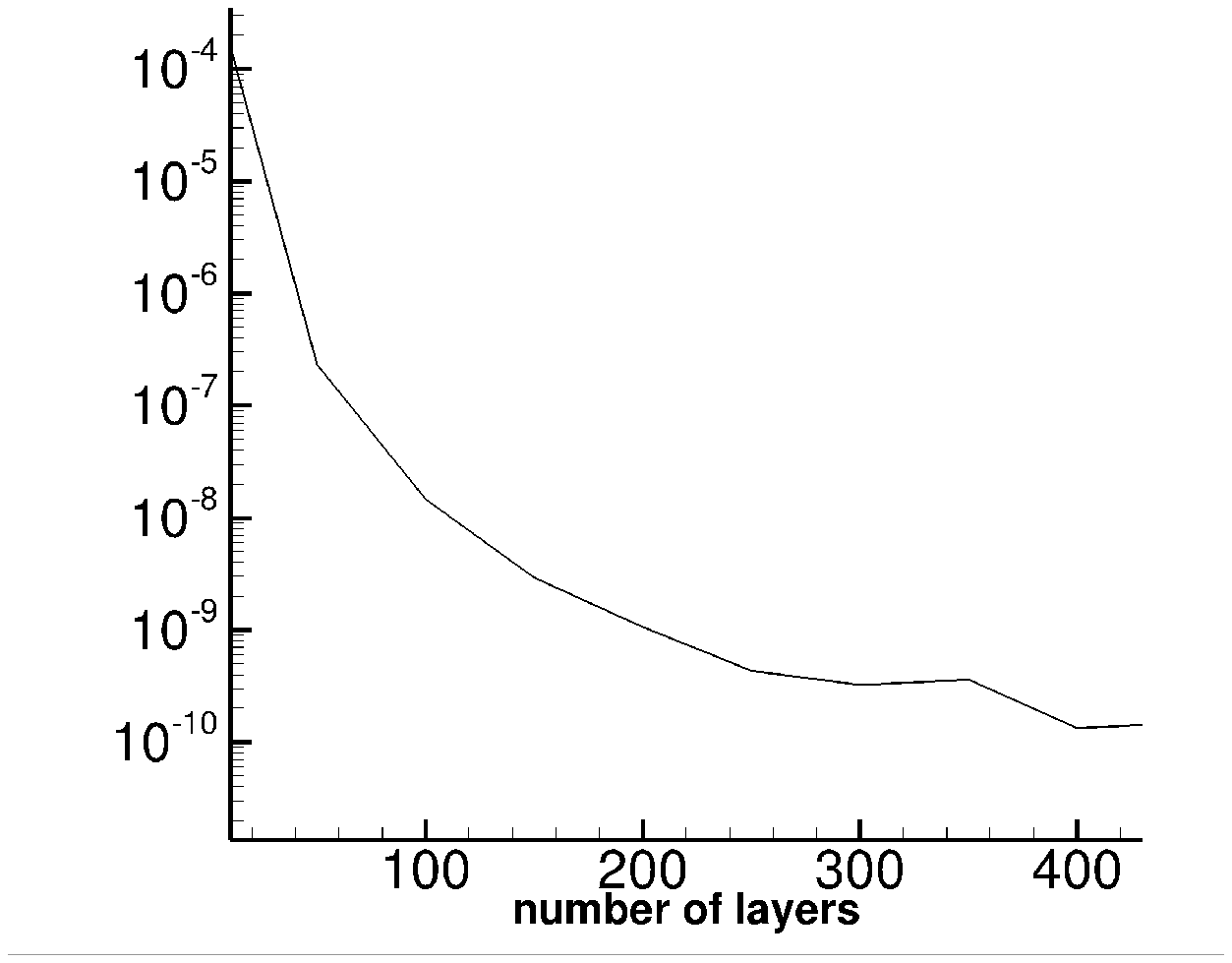}}
	\caption{Convergence of $e_{\psi}$ as the increase of length of input data (Left) and number of layers (Right). For the left one, we fix the number of layers to be 200 and for the right one, the number of neurons in each hidden layer is  chosen to be 20.}
	\label{fig:conver_exp}
\end{figure}

In this example, we consider \eqref{model_periodic_bc} with $\beta=-1$ , $\gamma=1$, and $\Omega=[-10,10]$. Moreover, the true potential has the exponential form as follows,
$$
V(x)=4x^2-\exp(-2x^2).
$$
Then the exact solution $\psi(x,t)$ is has the following form:
$$
\psi(x, t)=\exp(-x^2+2it).
$$
We first show the convergence of $e_{\psi}$ as the number of collocation points ($M$) and layers ($N$) in Figure \ref{fig:conver_exp}.  As for the left one, we fix the number of layers $N$ to be 200 and increase the number of collocation points $M$ in data.  One can see that the by increasing $M$ to 20, the error first shaply decrease to about $10^{-9}$ and become stable. On the other hand, with $M$ fixed as $200$, one can observe a steady decay in $e_{\psi}$ as using more layers, equivalently, smaller $\Delta t$. Hence the performance of the NLS-Net is consistent with the Theorem \ref{thm: conv_nn}. Since in our case, we deal with the inverse problem by iteratively solving the forward process, the effect of NLS-Net in approximating exact solution $\psi$ serve as a basis in finding the true potential.

\begin{figure}[htbp!]
	\centering
	\subfigure{\includegraphics[width=0.45\textwidth]{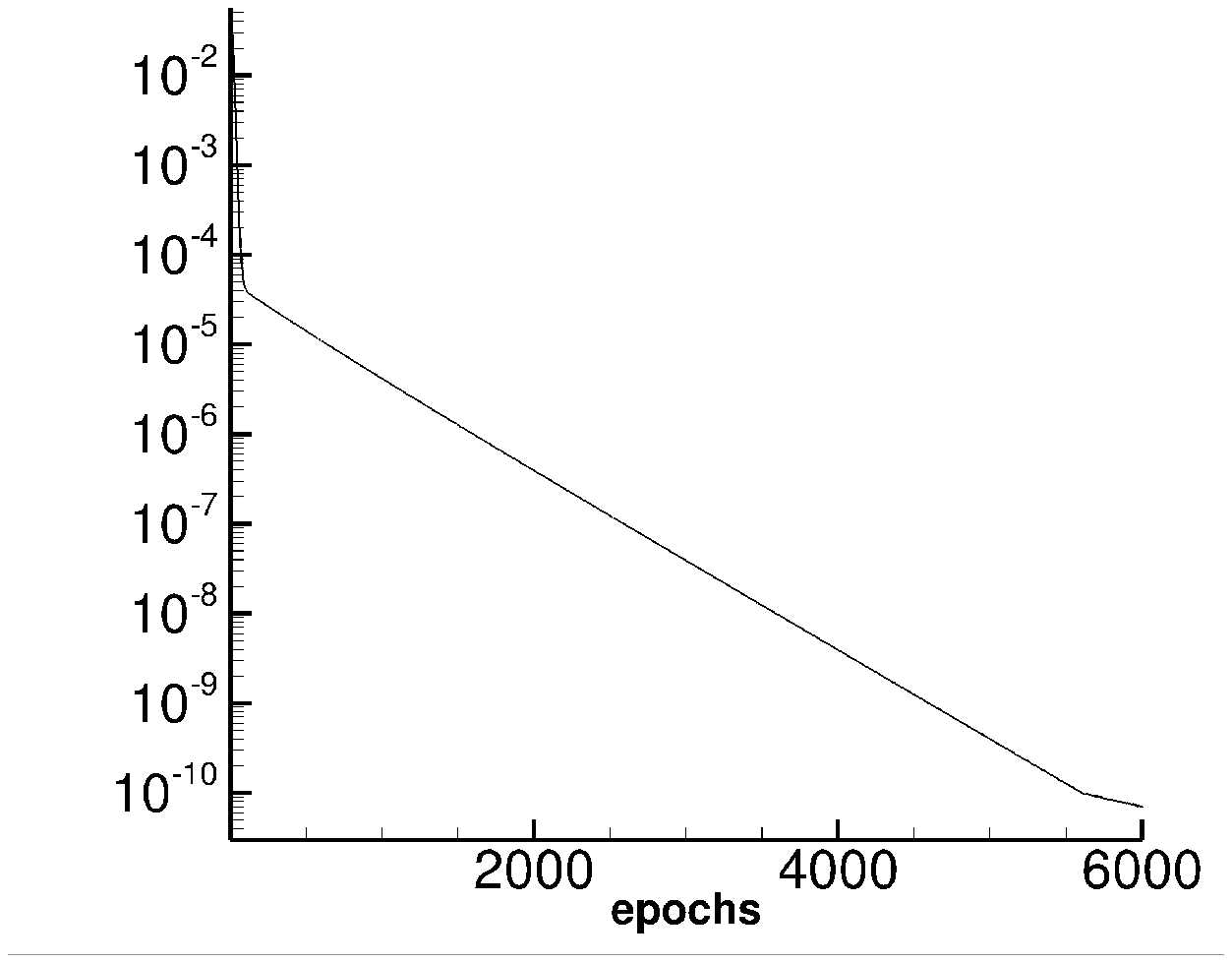}}
		\subfigure{\includegraphics[width=0.45\textwidth]{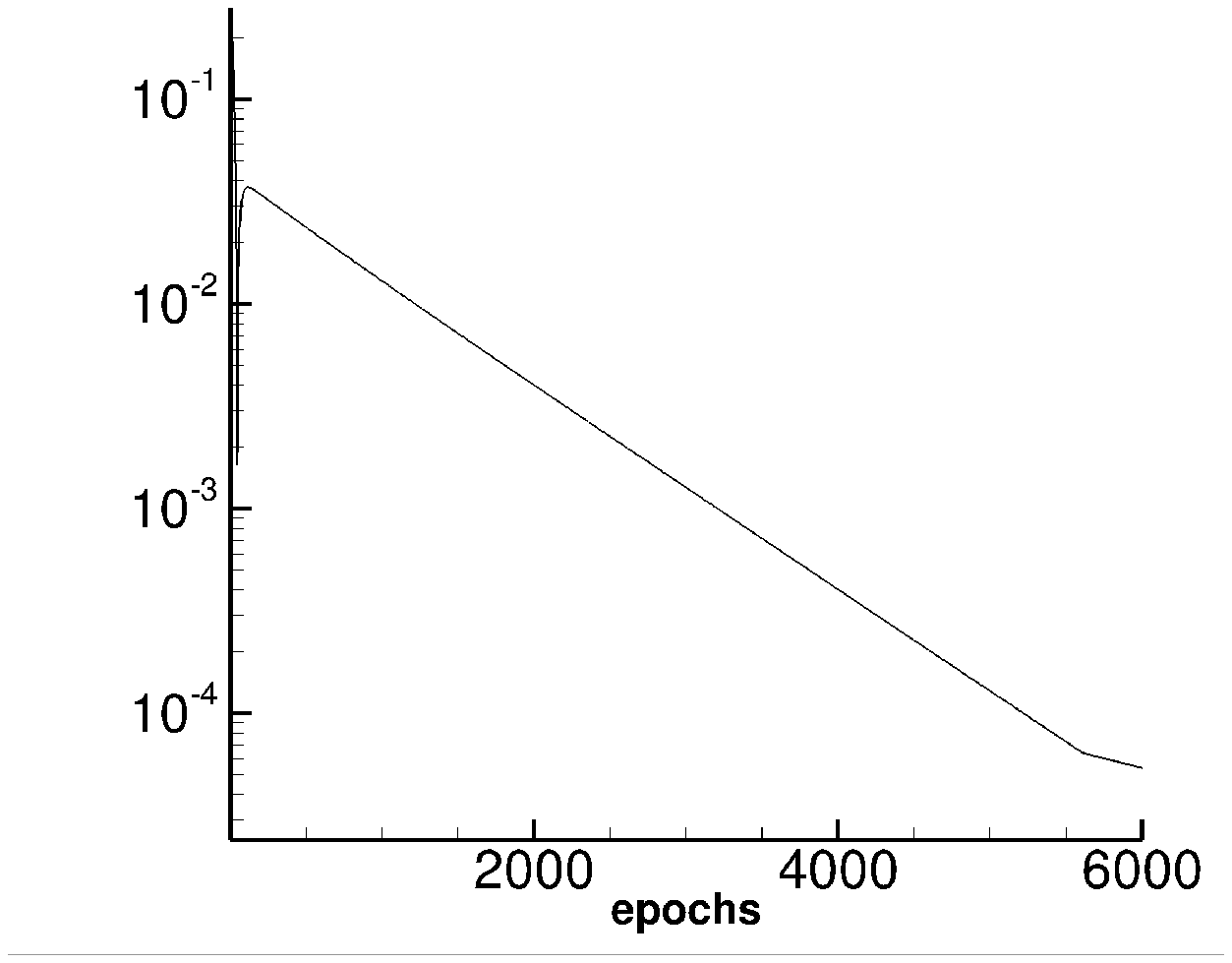}}
	\caption{Dynamics of loss function $e_{\psi}$ (Left) and potential error $e_{V}$ (Right) during training in \textbf{Example 1}. Here, the regularization coefficient $\lambda$ is zero, i.e., $J_{\bc}$ is the same as $e_{\psi}$. Gradient descent is used with initial learning rate $1.5$, which is multiplied by 0.99 every 2000 epochs. When $e_{\psi}$ decreases to be lower than $10^{-10}$, learning rate is multipied by $0.4$ each step.  }
	\label{fig:exp_error}
\end{figure}

We further present the training performance in Figure \ref{fig:exp_error}. The left subfigure describes the dynamics of $e_{\psi}$ in the training process while the right one plots $e_{V}$. Here we set the regularization coefficient $\lambda$ to be zero, which means $J_{\bc}=e_{\psi}$. Hence, we simply use gradient descient instead of proximal gradient descient. In particular, the initial learning rate is 1.5 and for each 2000 epochs, it decays by multiplying 0.99. Besides, we set a tolerence $\tau$ to be $10^{-10}$ as a creterion to further decrease the learning rate. Once $e_{\psi}$ is below this tolerance, we multiple learning rate by 0.4 for each remaining step. For the network, it is composed of four layers and with input $\tilde{\Psi}_M^0$ containing 4000 uniform collocation points in $[-10,10]$, i.e, $M=4000$.  Both the dynamics of $e_{\psi}$ and $e_{V}$ are plotted with logarithmic scale in y axis, hence one can observe a linear decay.  We remark that a change in slope of lines near the end of training is resulted by exponential decreasing the learning rate since $e_{\psi}$ already drops below a given small tolerence.

\begin{figure}[htbp!]
	\centering
	\subfigure{\includegraphics[width=0.45\textwidth]{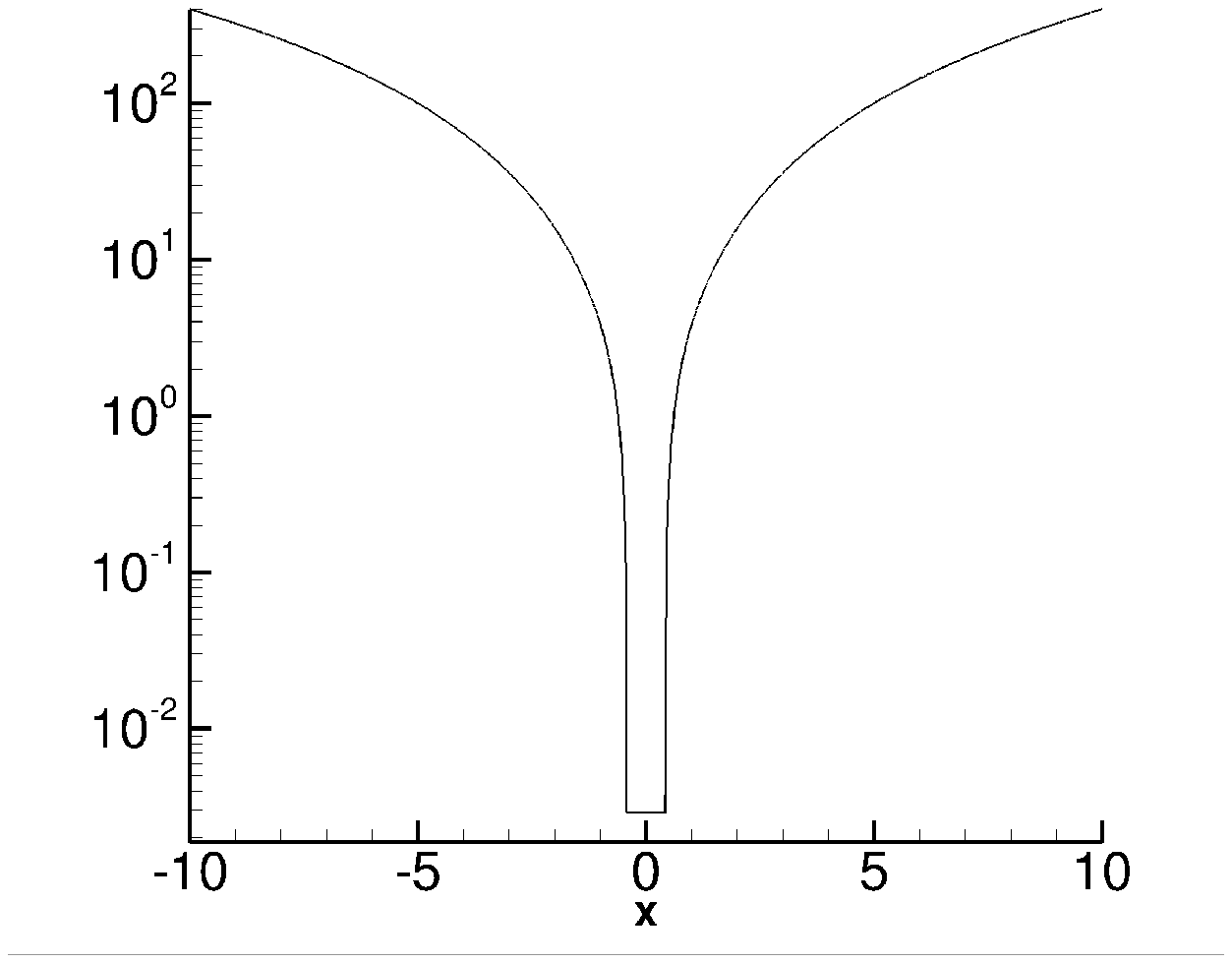}}
	\subfigure{\includegraphics[width=0.45\textwidth]{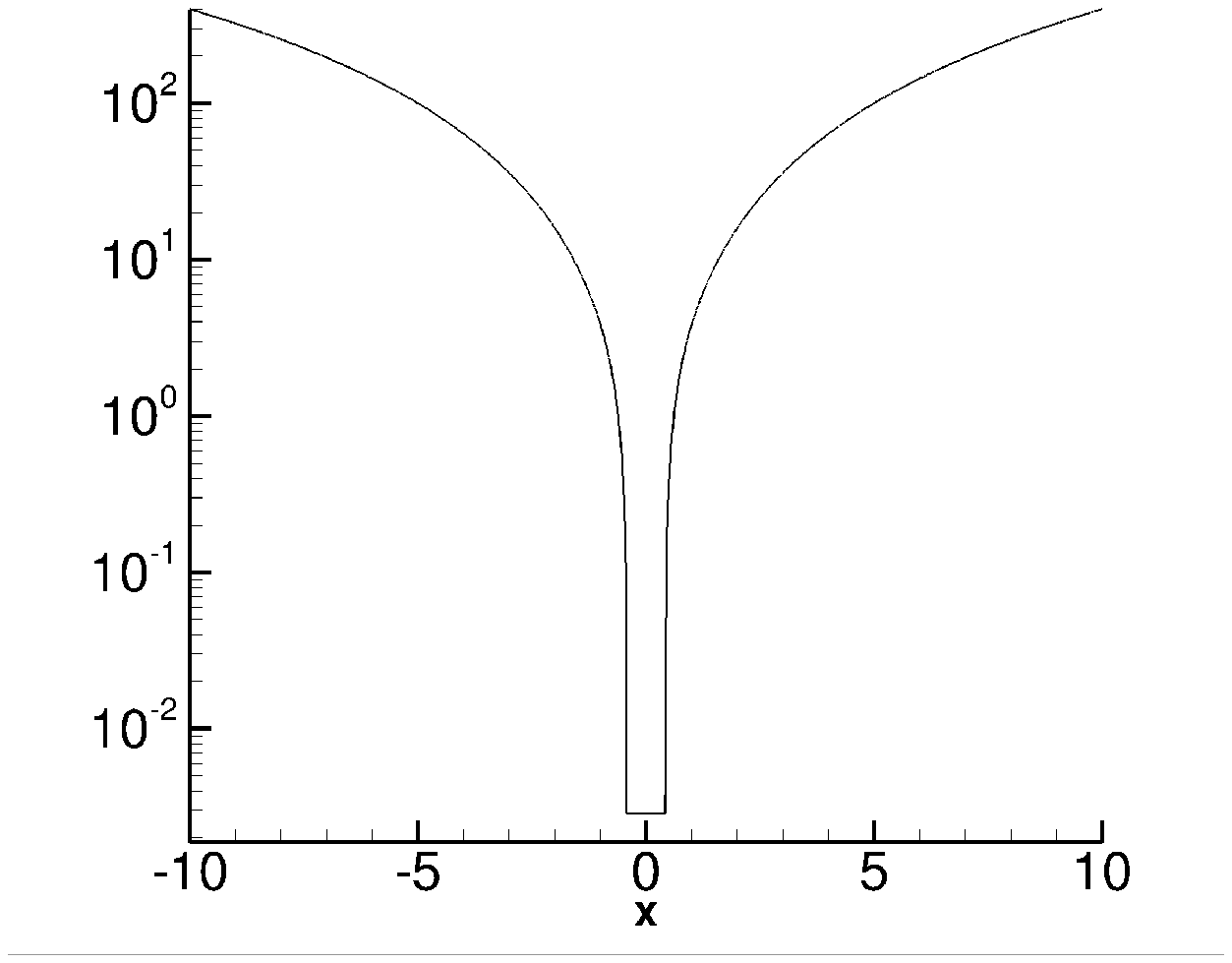}}
	\caption{A comparison of numerical (Left)and exact (Right) potential in Example 1 after training 6000 epochs. The spatial domain $[-10,10]$ is descritized by 4000 uniform collocation points and the time domain $[0,T]$ is partitioned by $\Delta t = T/4$. }
	\label{fig:exp_potential com}
\end{figure}

Finally we compare the numerical (Left) and exact (Right) potential in Figure \ref{fig:exp_potential com}. The domain $[-10,10]$ is discretized by 4000 uniform collocation points and we use 6000 epochs to training the potential. After that, one can observe that the apprximation of potential  is almost impossible to be distinguished from the goundtruth. 
\\
\newline
\textbf{Example 2.} Gross-Pitaevskii (GP) equation.

\begin{figure}[htbp!]
	\centering
	\subfigure{\includegraphics[width=0.45\textwidth]{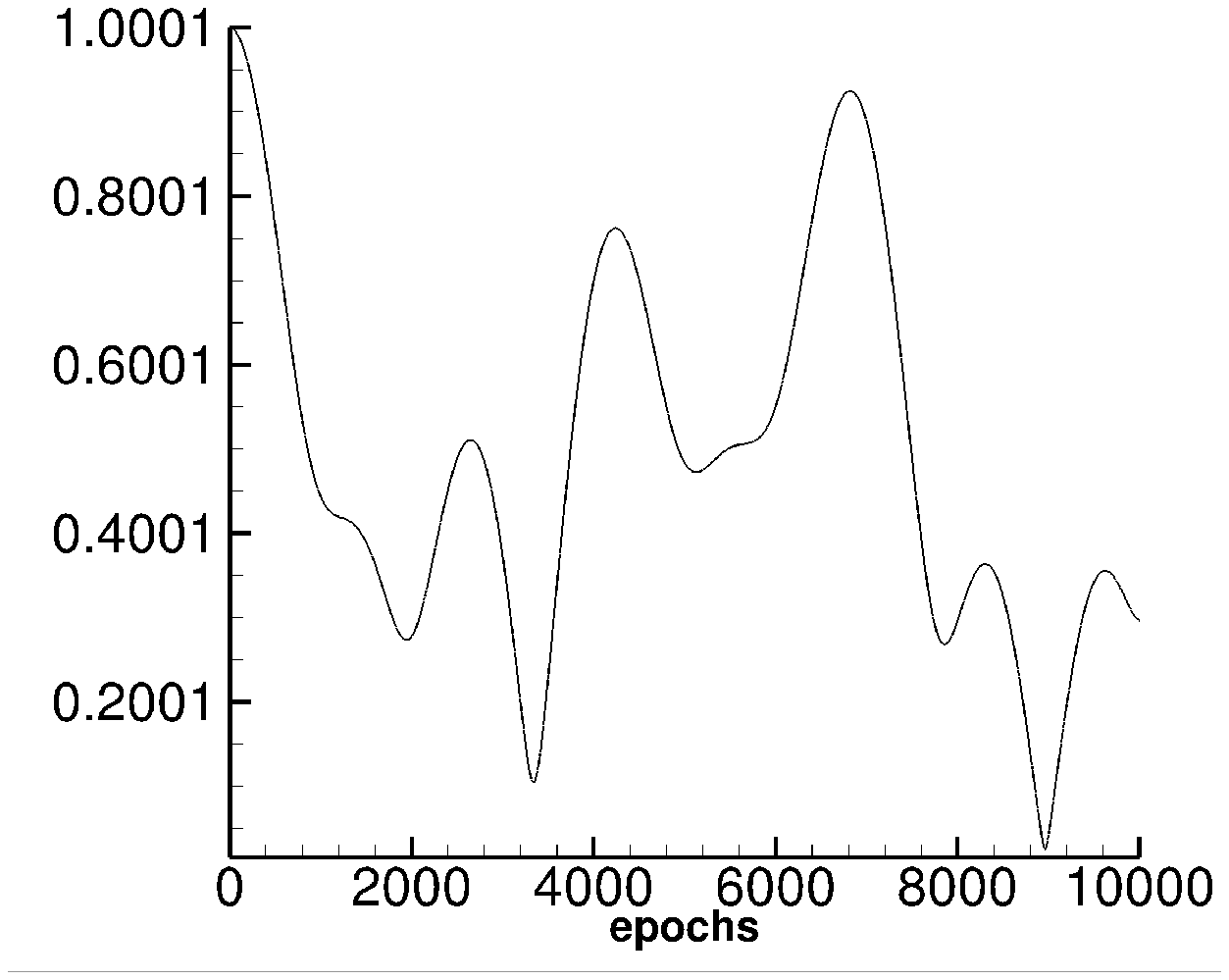}}
	\subfigure{\includegraphics[width=0.45\textwidth]{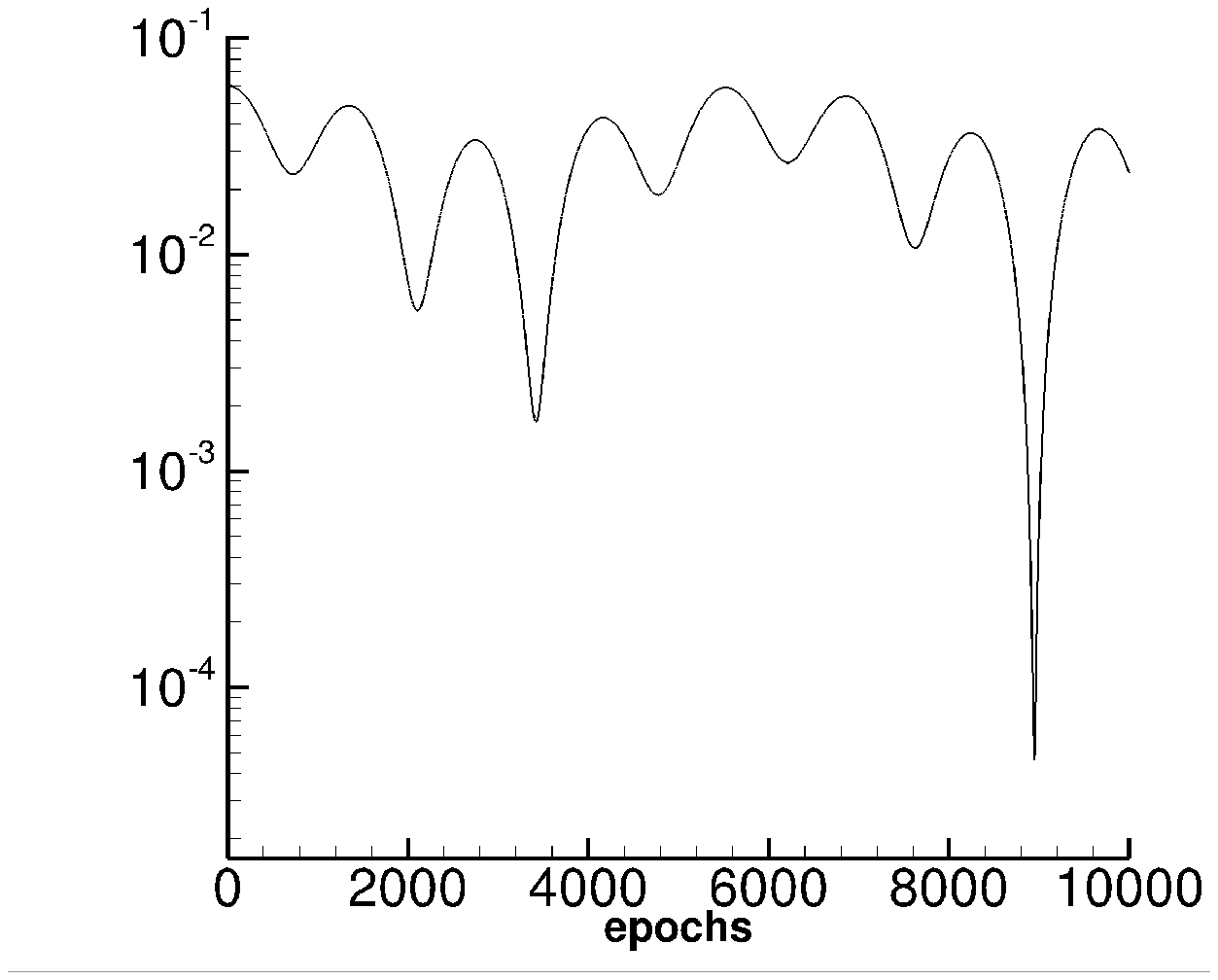}}
	\caption{Dynamics of loss function $e_{\psi}$ (Left) and potential error $e_{V}$ (Right) during training in Example 2. Here, the regularization coefficient $\lambda$ is zero, i.e., $J_{\bc}$ is the same as $e_{\psi}$. Gradient descent is used with initial learning rate $5*10^{-6}$, which is multiplied by 0.9 every 3000 epochs.  }
	\label{fig:cos_error_noreg}
\end{figure}

\begin{figure}[htbp!]
	\centering
		\subfigure{\includegraphics[width=0.45\textwidth]{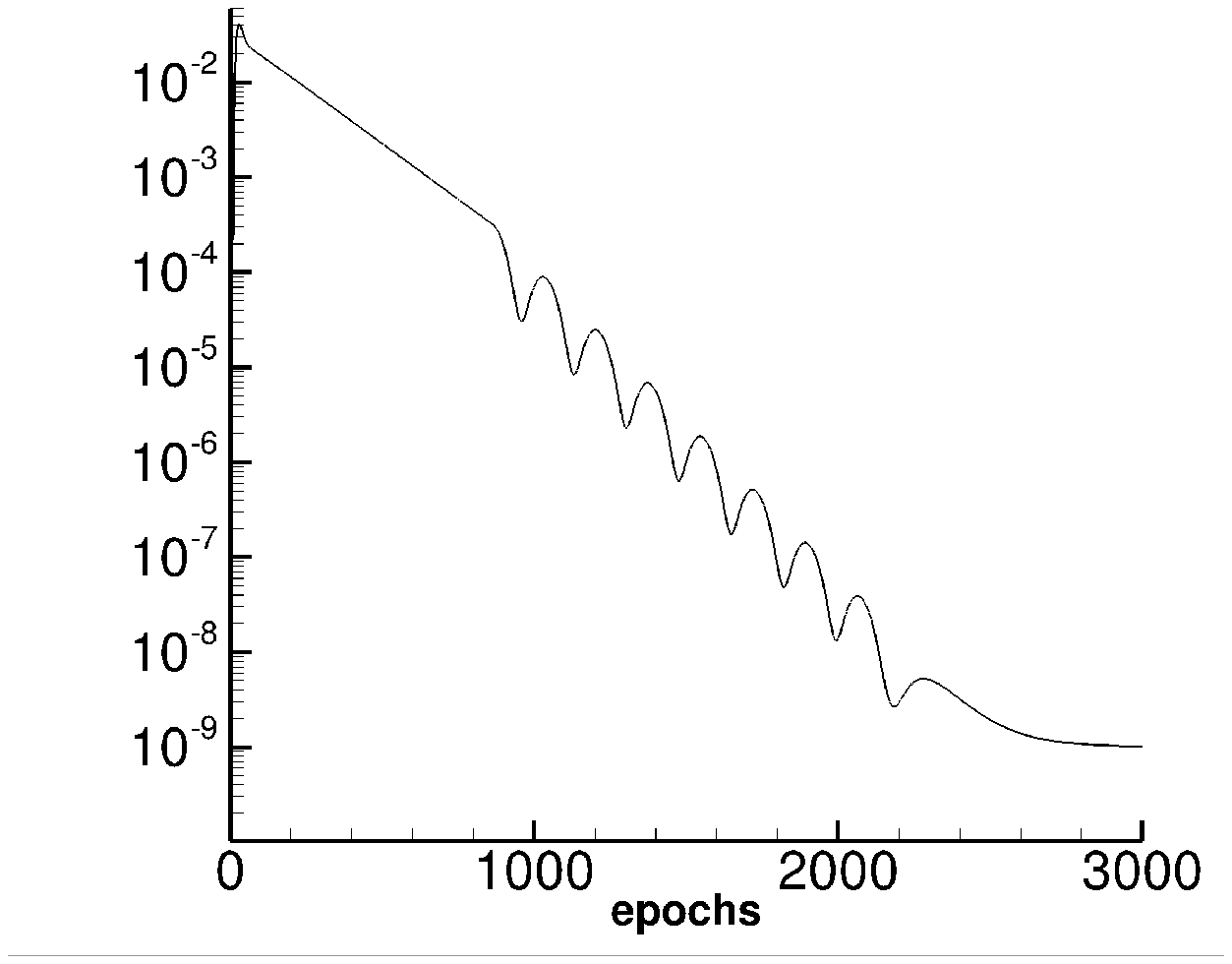}}
	\subfigure{\includegraphics[width=0.45\textwidth]{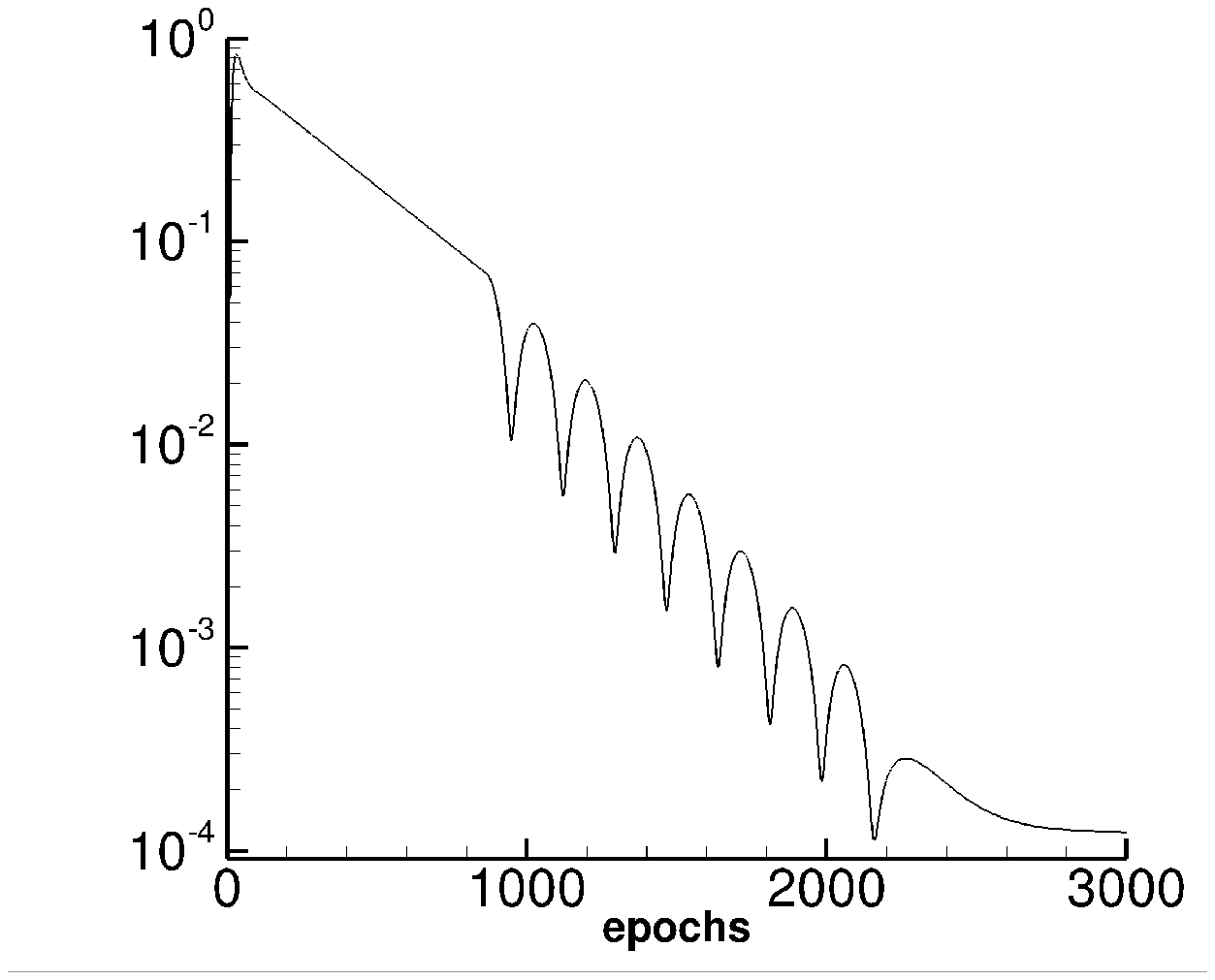}}
	\caption{Dynamics of loss function $e_{\psi}$ (Left) and potential error $e_{V}$ (Right) during training in \textbf{Example 2}. Here, the regularization coefficient $\lambda$ is 20. Proximal gradient descent is used with initial learning rate $5*10^{-3}$, which is multiplied by 0.9 every 3000 epochs.  When $e_{\psi}$ decreases to be lower than $10^{-10}$, learning rate is multiplied by $0.95$ each step.  }
	\label{fig:cos_error}
\end{figure}

In this example, we consider the GP equation, which is special case of \eqref{model_periodic_bc} with periodic potential. In particular, we set $\beta=\frac{1}{2}, \gamma=-1$,  the spatial domain $\Omega=$and $V(x)=-\cos ^2(x)$. For this equation, the exact solution \cite{wang2005numerical} is 
$$
\psi(x, t)=\sin (x) \exp (-\frac{3ti}{2}).
$$
First of all, we show the dynamics of $e_{\psi}$ and $e_{V}$ in the training process. To obtain a relatively good approximation, we choose $\lambda=20$ in \eqref{loss}. As a comparison, we also present the case without regularization term in Figure \ref{fig:cos_error_noreg}. In this case, we simply use SGD. Even with a very small learning rate $5*10^{-6}$, the loss $e_{\psi}$ still oscillates, which may imply a relatively high level of nonconvexity in $e_{\psi}$. Corresponding, the error $e_{V}$ also presents a similar trend during training. Both of them shows that it is relatively difficult to attain the optimal solution with the only term $e_{\psi}$ in $J_{\bc}$.  On the other hand, with a nonzero regulatization term, we can obtain a near exponential decay in both $e_{\psi}$ and $e_V$, which is shown in Figure \ref{fig:cos_error}. To obtain this result, we use proximal gradient descent. In particular, we use SGD to tackle the $l^2$ term and a proximal operator is utilized to deal with the $l^1$ norm, i.e. the second term in $J_{\bc}$. For SGD, we begin with a learning rate $5*10^{-3}$ and it decays by rate 0.1 every 3000 epochs. Similarly to the implementation in the last example, we shrink the learning by multiplying it with 0.95 each step after $e_{\psi}$ has satisfied the condition that $e_{\psi}<10^{-10}$. We remark that even though there are still regular oscillations in the dynamics because of SGD, both $e_{\psi}$ and $e_V$ display a steady decrease, which shows the effect of regularization term in $J_{\bc}$.

\begin{figure}[htbp!]
	\centering
	\subfigure{\includegraphics[width=0.45\textwidth]{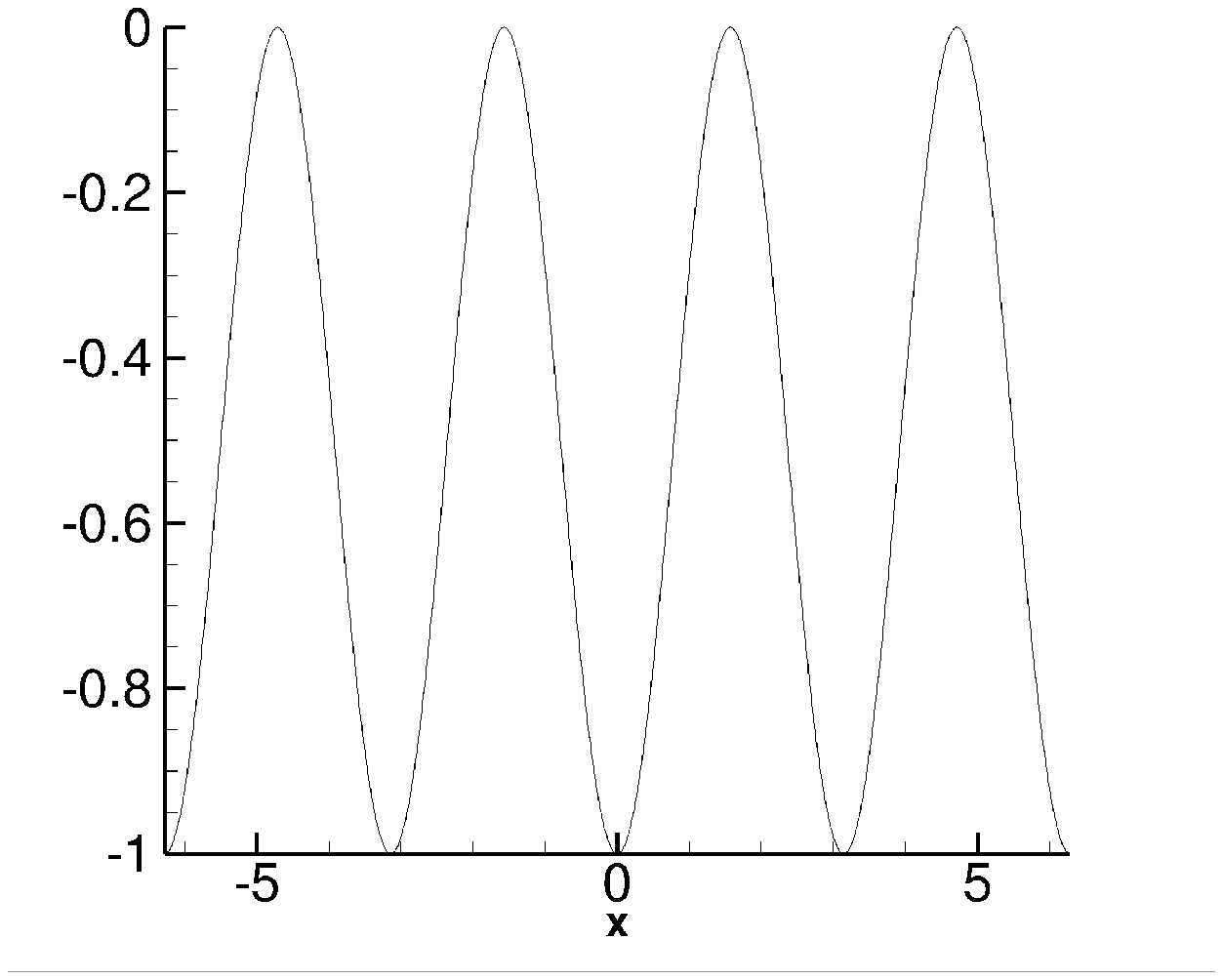}}
	\subfigure{\includegraphics[width=0.45\textwidth]{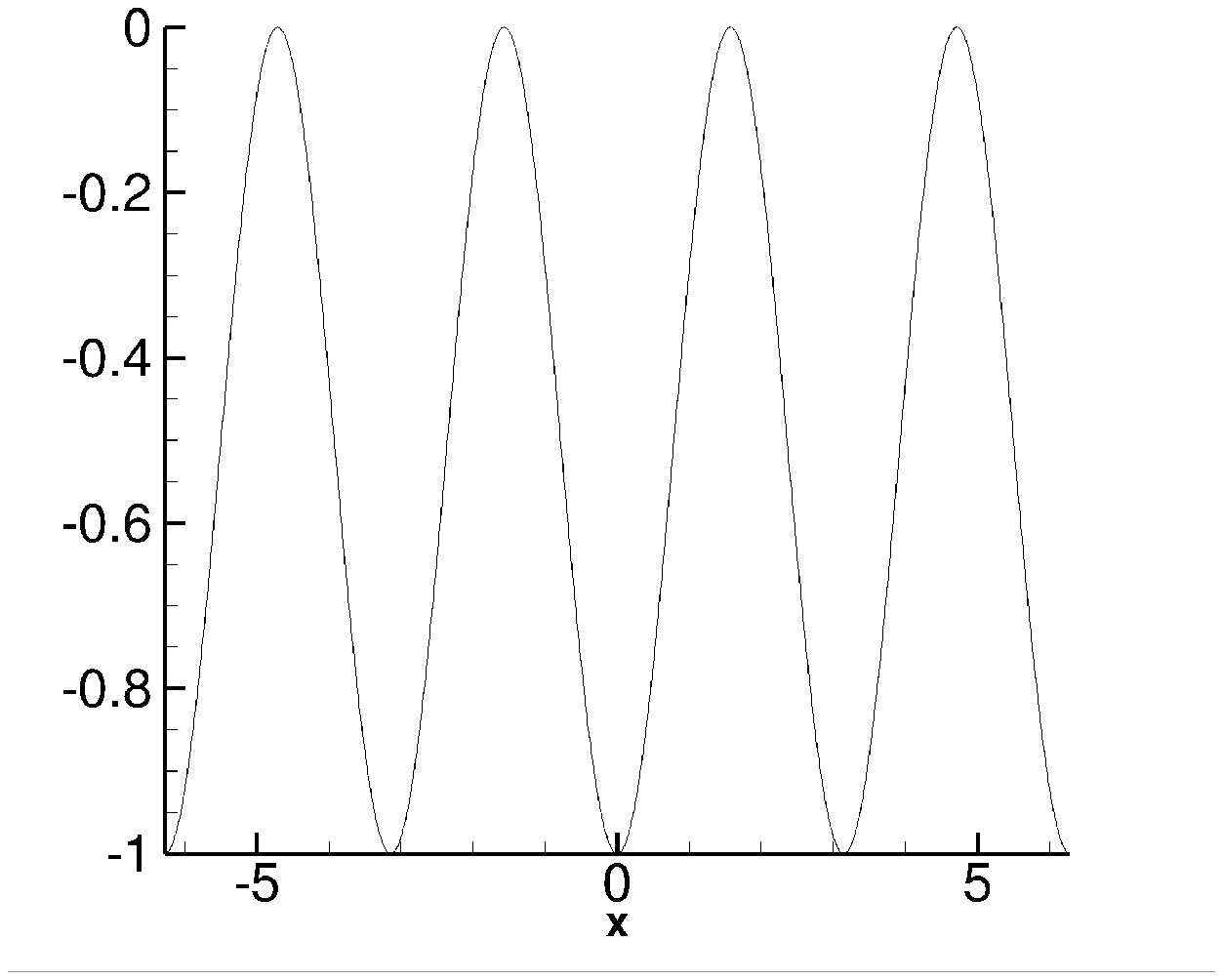}}
	\caption{A comparison of numerical (Left)and exact (Right) potential in Example 2 at final time $T=1$. The spatial domain $[-10,10]$ is descritized by 4000 uniform collocation points and the time domain $[0,T]$ is partitioned by $\Delta t = T$. }
	\label{fig:cos_potential com}
\end{figure}

In Figure \ref{fig:cos_potential com}, we compare the approximation of potential which is plotted on the left with the groundtruth on the other side. One can easily see that the numerical potential is approximately identical to the reference, which shows that the training process is effective.
\\
\newline
\textbf{Example 3.} A coupled system of NLSE.

\begin{figure}[htbp!]
	\centering
	\subfigure{\includegraphics[width=0.5\textwidth]{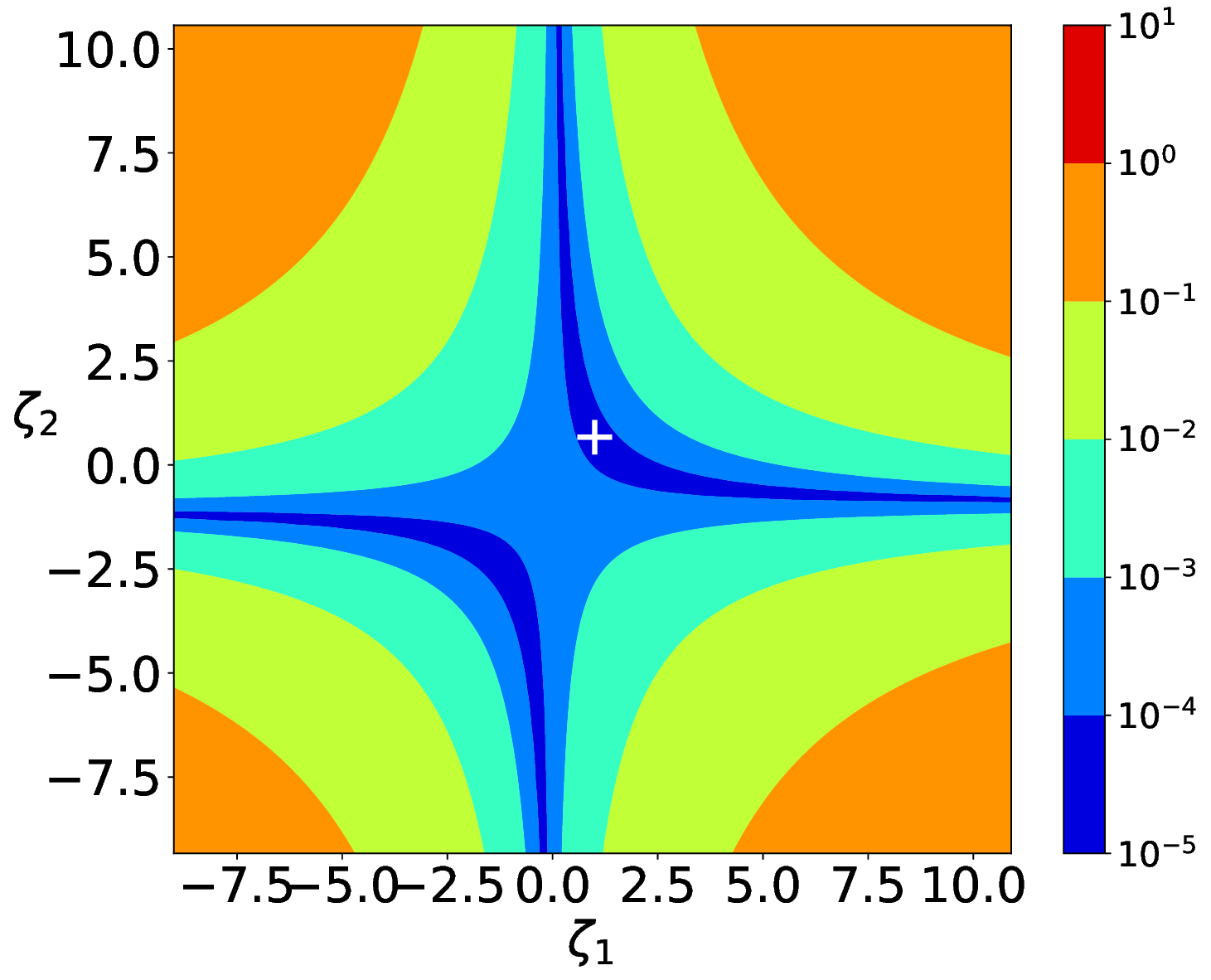}}
	\label{fig:landscape_couple}
\end{figure}

In this example, we explore a coupled nonlinear Schr\"{o}dinger equations as follows. For $x\in[a,b]$, $t\in[0,T]$,
\begin{eqnarray*}
\begin{aligned}
	& \mathrm{i} \frac{\partial \psi_j}{\partial t}=-\frac{1}{2} \frac{\partial^2 \psi_j}{\partial x^2}+V_j(x) \psi_j+i \alpha_j \frac{\partial \psi_j}{\partial x}+f_j\left(\left|\psi_1\right|^2,\left|\psi_2\right|^2\right) \psi_j, \quad j=1,2, \\
	& \psi_j(x, 0)=\psi_j^0(x), \quad j=1,2 \\
	& \psi_j(a, t)=\psi_j(b, t)=0, \quad j=1,2 ,
\end{aligned}
\end{eqnarray*}
where $f_1(x, y)=x+\zeta_1 y$, $f_2(x, y)=\zeta_2 x+y.$ Here, $\zeta_1$ and $\zeta_2$ are scalars.  We can know that it admits one exact solution \cite{wang2005numerical} as follows when $ \zeta_1= \zeta_2$, $V_1(x)=V_2(x)=0$, and $-\alpha_1=\alpha_2=0.5$. The overall loss is defined as $J_{\bc}= e_{\psi_1}+e_{\psi_2}$, where $\bc=[\zeta_1,\zeta_2]$.
\begin{eqnarray}
\begin{aligned}
	& \psi_1^{\text {exact }}(x, t)=\sqrt{\frac{2 \alpha}{1+\zeta}} \operatorname{sech}(\sqrt{2 \alpha}(x-v t)) \mathrm{e}^{\mathrm{i}\left[(v-\delta) x-\left(\frac{2^2-\delta^2}{2}-\alpha\right) t\right]} \\
	& \psi_2^{\text {exact }}(x, t)=\sqrt{\frac{2 \alpha}{1+\zeta}} \operatorname{sech}(\sqrt{2 \alpha}(x-v t)) \mathrm{e}^{\mathrm{i}\left[(v+\delta) x-\left(\frac{2^2-\delta^2}{2}-\alpha\right) t\right]}
	\label{sol_couple}
\end{aligned}
\end{eqnarray}
In our numerical simulation, we use the following setting:
\begin{align*}
a=-20, \quad b=80, \quad v=0., \quad \delta=0.5, \quad \alpha=1.0, \quad \zeta_1= \zeta_2=\frac{2}{3}.
\end{align*}
The difference of this example is that the potential is a function of forward solutions $\psi_1$ and $\psi_2$, while in the first two examples, we have the prior information that the unknown potential is only a function of $x$. However, as a compensation, the form of potential is specified and we  only need to determine $\zeta_1$ and $\zeta_2$ rather than find a unknown function. We remark that even though the reference solutions defined in \eqref{sol_couple} correponds to the special case $\zeta_1= \zeta_2$, it is not reasonable to assume this information hence we need to determine both $\zeta_1$ and $\zeta_2$. However, since there are still only two unknown scalars,   we can easily plot the landscape of loss function $J_{\bc}$ as a guide to seek the optimal solution. In Figure \ref{fig:landscape_couple}, we show the contour of loss with respect to $\zeta_1$ and $\zeta_2$. First of all, it presents a reflection symmetry along the line $\zeta_1= \zeta_2$, which can offer the insight that the optimal $\zeta_1$ and $\zeta_2$ are almost equal. Besides, the loss function shows the poperty of local convexity. Hence once the $[\zeta_1,\zeta_2]$ arrives the basin of attraction, where we borrow the terminology in ecology, it is not demanding to attain the optimal solution. Both two observations can help us find the minimizer in a more efficient way.

Following the loss landscape, we further show the training performance. As before, we arrange the loss dynamics on the left and  errors of inverse problems on the other side. Since in this situation, we only need to determine the unknowns scalars, we use the absolute error defined as 
\begin{align*}
	e_{\zeta_j} = \frac{|\zeta_j^{\text{numerical}}-\zeta_j^{\text{exact}}|}{|\zeta_j^{\text{exact}}|}, \quad j=1,2.
	\end{align*}
Owing to the landscape analysis shown above, we can simply use gradient descient. Here we fix the learning rate to be 100. The reason for such a big learning rate may come from a small magnititude  of the gradient with respect to the parameters $\zeta_1$ and $\zeta_2$. However, in spite of the big learning rate, we still obtain an exponeital decay of loss as well as the concerned errors. One can observe a convergence before about 1000 epochs.
\begin{figure}[htbp!]
	\centering
		\subfigure{\includegraphics[width=0.45\textwidth]{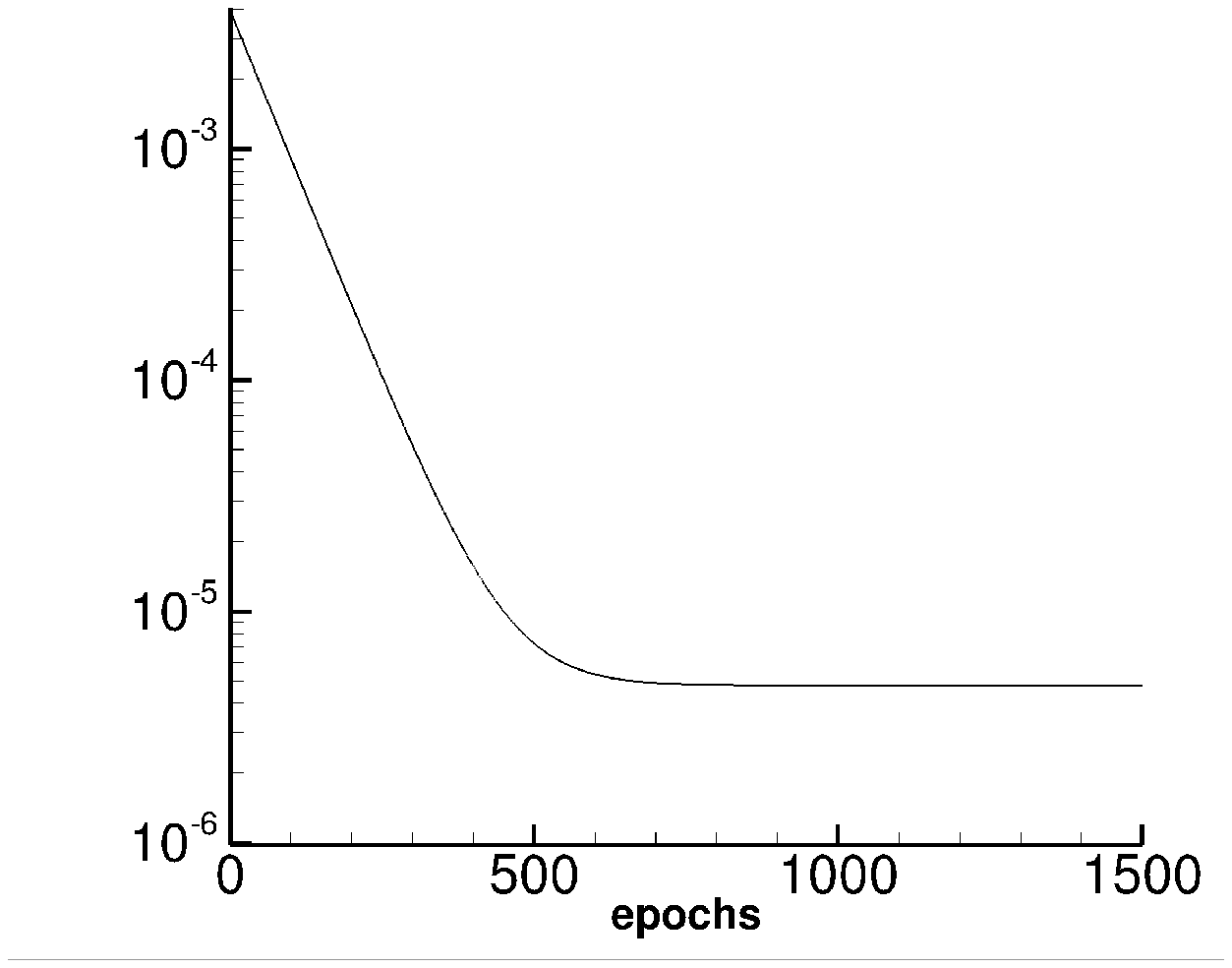}}
	\subfigure{\includegraphics[width=0.45\textwidth]{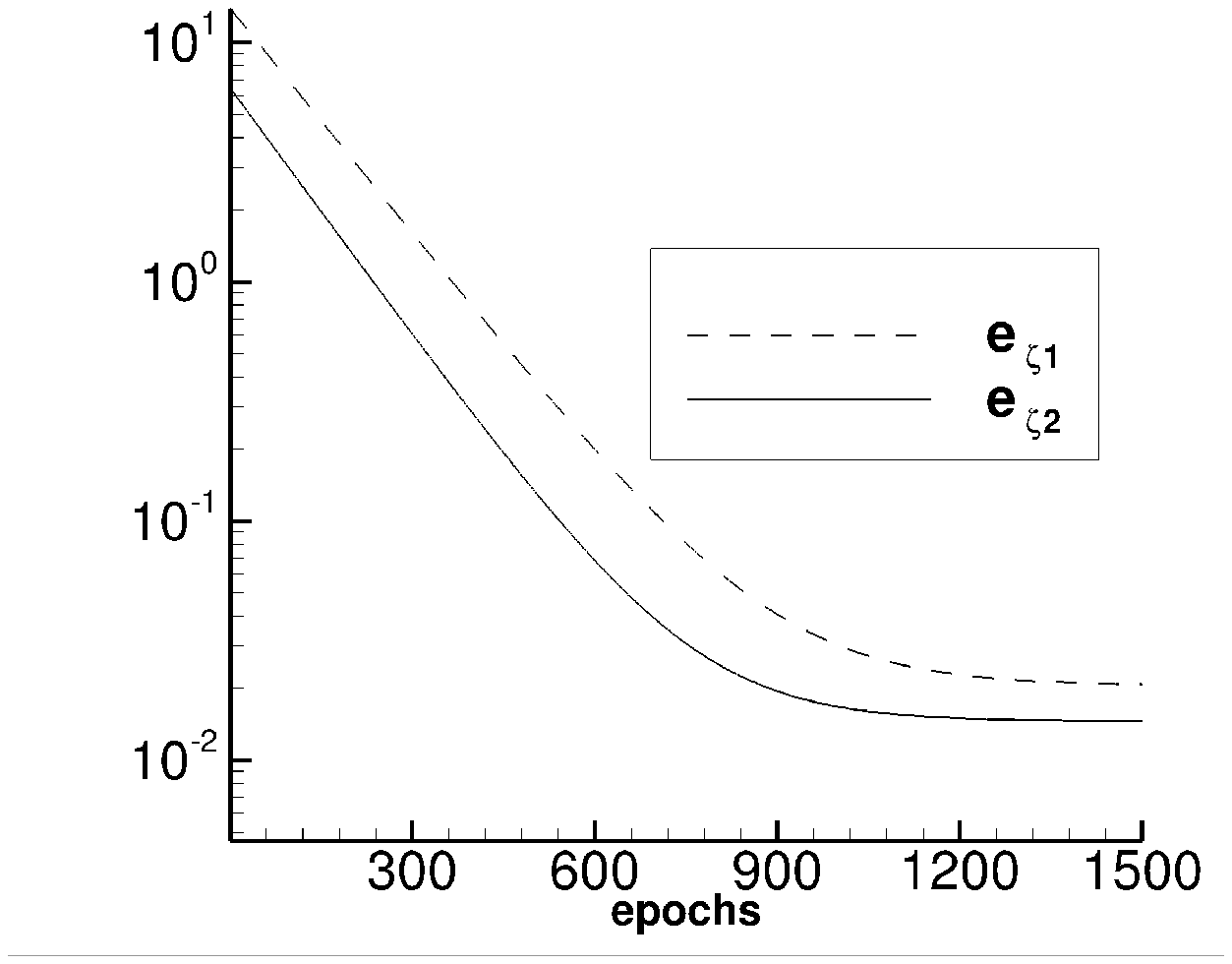}}
	\caption{Dynamics of loss function $e_{\psi}$ (Left) , errors of $\zeta_1$ and $\zeta_2$ (Right) during training in Example 3.  Gradient descent is used with fixed learning rate $100$. }
	\label{fig:couple_error}
\end{figure}	
\section{Conclusion remarks}
In this work, we solve some inverse problems of nonlinear Schr\"{o}dinger equations as a learning process of a special type of convolutional neural network (CNN). The first novelty of this work comes from  the special architecture of the neural network, which is  designed based on the equation and the operator splitting method. In particular, the convolution layers serve as a linear operator while the activation functions constitude the nonlinear part. A composition of these two  forms the time-marching process.  Owing to this design, layers and activation functions have physical interpretations.  An intermediate convenience is that one can easily choose the type of activation functions from the information of equations. Besides, the existence of an optimal network with this special structure can be guaranteed since the neural network can be directly derived from the equation. We remark that the existence is not a direct result of the famous universal approximation theorem for neural networks since the proposed network has a special architecture while the theorem does not specify the structure of the optimal network. Moreover, based on the existing theoretical result of operator splitting, one can easily guarantee the convergence of the neural network solution with respect to the depth and width. Overall, the particular structure derived from the equation enhances the explainability of network. Instead of using neural network as a black box, the inherited physics sheds light on the approximation effect of the neural network.

For the inverse problem, we are interested in approximating the potential function in the Schr\"{o}dinger equations. Different from the coefficient inverse problem, function approximation is more complicated. For example, one needs a neural network to approximate a function while for coefficients, some neurons are sufficient. To this end, we use  a library-search method to transform the original problem to a coefficient inverse problem, which is the second contribution. In particular, we construct a library that contains a set of basis functions then we can use a linear combination of the bases to express the unknown potential. In other words, we project the solution of inverse problems to an approximation space so we only need to determine the coefficients with respect to the basis functions. In this way, the complexity of neural network can be effectively reduced. However, to approximate a function space, we need to incorporate sufficiently many bases to obtain a relatively good accuracy. On the other hand, the form of potential is not so complicated, which may be a linear combination of a small number of different functions. Hence, the coefficient vector correponding to basis functions should have some level of sparsity, which motivates the need of using compressed sensing. 
	\bibliographystyle{plain}
\bibliography{references}
\end{document}